\documentclass{cocv}

\usepackage{amsfonts}
\usepackage{graphicx}
\usepackage{epstopdf}
\ifpdf
  \DeclareGraphicsExtensions{.eps,.pdf,.png,.jpg}
\else
  \DeclareGraphicsExtensions{.eps}
\fi

\usepackage{amssymb}
\usepackage{amsfonts}
\usepackage{bbm}
\usepackage{xcolor}

\usepackage{xspace}

\usepackage[font={small,it}]{caption}
\usepackage{algorithm}
\usepackage[noend]{algpseudocode}
\usepackage{hyperref}
\usepackage{enumitem}

\newcommand{\firstRev}[1]{\textcolor{black}{#1}}

\begin{document}
\title{Sequential Convex Programming for \\
Non-Linear Stochastic Optimal Control}

\author{Riccardo Bonalli}\address{The author is with the Université Paris-Saclay, CNRS, CentraleSupélec, Laboratoire des signaux et systèmes, 91190, Gif-sur-Yvette, France. \textit{Email:} {\tt riccardo.bonalli@centralecupelec.fr} and {\tt rbonalli@stanford.edu}.}
\author{Thomas Lew$^2$,}
\author{Marco Pavone}\address{The authors are with the Department of Aeronautics and Astronautics, Stanford University, Stanford, CA 94305. \textit{Emails:} {\{\tt thomas.lew, pavone\}@stanford.edu}. 
This research was supported by the National Science Foundation under the CPS program (grant \#1931815).}

\begin{abstract}
  This work introduces a sequential convex programming framework \firstRev{for} non-linear, finite-dimensional stochastic optimal control, where uncertainties are modeled by a multidimensional Wiener process. We prove that \firstRev{any accumulation point of the sequence of iterates generated by sequential convex programming is} a candidate locally-optimal solution for the original problem in the sense of the stochastic Pontryagin Maximum Principle. \firstRev{Moreover, we provide sufficient conditions for the existence of at least one such accumulation point.} We then leverage these properties to design a practical numerical method for solving non-linear stochastic optimal control problems based on a deterministic transcription of stochastic sequential convex programming.
\end{abstract}

\begin{resume}
    Ce travail introduit un cadre de programmation s\'equentielle convexe \firstRev{pour le contr\^ole optimal de syst\`emes stochastiques non-lin\'eaires de dimensions finies.} Nous prouvons que \firstRev{tout point d'accumulation de la suite des it\'er\'ees g\'en\'er\'ee par la programmation séquentielle convexe est} un candidat \`a une solution localement optimale du probl\`eme originel, au sens du Principe du Maximum de Pontriaguine stochastique. \firstRev{De plus, nous d\'eveloppons des conditions suffisantes pour l'existence d'au moins un de ces points d'accumulation.} Nous exploitons ces propri\'et\'es afin de concevoir une m\'ethode numérique pratique r\'esolvant des probl\`emes de contr\^ole stochastique non-lin\'eaires par une reformulation d\'eterministe de la programmation s\'equentielle convexe.
\end{resume}



\maketitle

\section{Introduction}

Over the past few decades, the applied control community has devoted increasing attention toward the optimal control of stochastic systems. The general formulation of this problem entails steering a dynamical system from an initial configuration to a final configuration while optimizing some prescribed performance criterion (e.g., minimizing some given cost) and satisfying constraints. This dynamical system is also subject to uncertainties which are modeled by Wiener processes and may come from unmodeled and/or unpredicted behaviors, as well as from measurement errors. Active development of new theoretical and numerical methods to address this problem continues, and the subject already has a rich literature.

We can %
classify the existing works into two main categories.

The first category consists of contributions that focus on Linear Convex Problems (LCPs), i.e., whose dynamics are linear and whose costs are convex in both state and control variables. An important class of LCPs is given by Linear Quadratic Problems (LQPs) in which costs are quadratic in both state and control variables and for which the analysis of optimal solutions may be reduced to the study of an %
algebraic relation known as Stochastic Riccati Equation (SRE) \cite{potter1965,bismut1976,peng1992,tang2003}. Efficient algorithmic frameworks have been devised to numerically solve LCPs, ranging from local search \cite{kleindorfer1973} and dual representations \cite{rockafellar1990,kuhn2011}, to deterministic-equivalent reformulations \cite{bes1989,berret2020}, among others. In the special case of LQPs, those techniques may be further improved by combining SRE theory with semidefinite programming \cite{rami2000,yao2001}, finite-dimensional approximation \cite{Bertsimas2007,damm2017}, or chaos expansion \cite{levajkovic2018}.

The second category of works deals with problems that do not enjoy any specific regularity, allowing non-linear dynamics or non-convex (therefore non-quadratic) costs. Throughout this paper, we call these Non-Linear Problems (NLPs). It is unquestionable that NLPs have so far received less attention from the community than LCPs, especially since the analysis of the former is usually more involved. Similar to the deterministic case, there are two main theoretical tools that have been developed to analyze NLPs: stochastic Dynamic Programming (DP) \cite{bellman1957,lions1983} and the stochastic Pontryagin Maximum Principle (PMP) \cite{pontryagin2018,kushner1964,Peng1990} (an extensive survey of generalizations of DP and the PMP may be found in \cite{Yong1999}). In the case of LQPs, one can show that DP and the PMP lead to SRE \cite{Yong1999}. DP provides optimal policies through the solution of a partial differential equation, whereas the necessary conditions for optimality offered by the PMP allow one to set up a two-point boundary value problem which returns candidate locally-optimal solutions when solved. Both methods only lead to analytical solutions in a few cases, and they can involve complex numerical challenges (the stochastic setting is even more problematic than the deterministic one, the latter being better understood for a wide range of problems, see, e.g., \cite{trelat2012,Bonalli2017Bis,bonalli2019_TAC}). This has fostered the investigation of more tractable approaches to solve NLPs such as Monte Carlo simulation \cite{shapiro2005,gobet2016}, Markov chain discretization \cite{kushner1990,kushner1991}, and deterministic (though non-equivalent) reformulations \cite{annunziato2013,berret2020}, among others. Importantly, many of the aforementioned approaches, e.g., \cite{kushner1990}, are often based on some sort of approximation of the original formulation and thus offer powerful alternatives to DP and the PMP, especially since they are more numerically tractable and can be shown to converge to policies satisfying DP or the stochastic PMP.

In this paper, our objective is to lay the theoretical foundations to leverage Sequential Convex Programming (SCP) for the purposes of computing candidate optimal solutions for a specific class of NLPs. SCP is among the most well-known and earliest approximation techniques for deterministic non-linear optimal control and, to the best of our knowledge, such an approach has not been extended to stochastic settings yet.
The simplest SCP scheme (which we consider in this work) consists of successively linearizing any non-linear terms in the dynamics and any non-convex functions in the cost and seeking a solution to the original formulation by solving a sequence of LCPs. %
This approach leads to two desirable properties, which jointly are instrumental to the design of efficient numerical schemes. First, one can rely on the many efficient techniques and associated software libraries that have been devised to solve LCPs (or LQPs depending on the shape of the original NLP). Second, as we will show in this paper, when this iterative process converges, it returns a strategy %
that satisfies the PMP related to the original NLP, i.e., a candidate optimum for the original formulation. %
Unlike existing methods such as \cite{kushner1990}, which introduce approximated formulations that are still non-linear, our approach offers the main advantage of requiring the solution to LCPs only. Specifically, we identify three key contributions:
\begin{enumerate}
    \item We introduce and analyze a new framework to compute candidate optimal solutions for finite-horizon, finite-dimensional non-linear stochastic optimal control problems. with control-affine dynamics and uncontrolled diffusion. This hinges on the basic principle of SCP, i.e., iteratively solving a sequence of LCPs that stem from successive linear approximations of the original problem.
    
    \item Through a meticulous study of the continuity of the stochastic Pontryagin cones of variations with respect to linearization, we prove that \firstRev{any accumulation point of the sequence of iterates generated by SCP is} a strategy satisfying the PMP related to the original formulation. In addition, by leveraging additional assumptions we prove that \firstRev{any such accumulation point always exists}, which in turn provides a ``weak'' guarantee of success for the method.
    
    \item Through an explicit example, we show how to leverage the properties offered by this framework to better understand what approximations may be adopted for the design of %
    efficient numerical schemes for NLPs\firstRev{, although we leave the theoretical investigation of the approximation error as future direction.}
\end{enumerate}

The paper is organized as follows. Section \ref{sec_problems} introduces notation and preliminary results and defines the stochastic optimal control problem of interest. In Section \ref{sec_SCP}, we introduce the framework of stochastic SCP and the stochastic PMP, and we state our main result of convergence. For the sake of clarity, in Section \ref{sec_proof} we retrace the proof of the stochastic PMP and introduce all the necessary technicalities we need to prove our main result of convergence, though we recall that the stochastic PMP is a well-established result and Section \ref{sec_proof} should not be understood as part of our main contribution. In Section \ref{sec_numerical}, we show how to leverage our analysis to design a practical numerical method to solve non-linear stochastic optimal control problems, and we provide numerical experiments. Finally, Section \ref{sec_conclusion} provides concluding remarks and perspectives on future directions.

\section{Stochastic Optimal Control Setting} \label{sec_problems}

Let $(\Omega,\mathcal{G},P)$ be a second--countable probability space and let $B_t = (B^1_t,\dots,B^d_t)$ be a $d$--dimensional Brownian motion with continuous sample paths starting at zero and whose filtration $\mathcal{F} \triangleq (\mathcal{F}_t)_{t \ge 0} = \big(\sigma(B_s : 0 \le s \le t)\big)_{t \ge 0}$ is complete. We consider processes that are defined within bounded time intervals. Hence, for every $n \in \mathbb{N}$, $\ell \ge 2$ and maximal time $T \in \mathbb{R}_+$, we introduce the space $L^{\ell}_{\mathcal{F}}([0,T] \times \Omega;\mathbb{R}^n)$ of progressive processes $x : [0,T] \times \Omega \rightarrow \mathbb{R}^n$ (with respect to the filtration $\mathcal{F}$) such that $\mathbb{E}\left[ \int^T_0 \| x(s,\omega) \|^{\ell} \; \mathrm{d}s \right] < \infty$, where $\| \cdot \|$ is the Euclidean norm. In this setting, for every $x \in L^{\ell}_{\mathcal{F}}([0,T] \times \Omega;\mathbb{R}^n)$ and $i=1,\dots,d$, the It\^o integral of $x$ with respect to $B^i$ is the continuous, bounded in $L^2$, and $n$--dimensional martingale in $[0,T]$ (with respect to the filtration $\mathcal{F}$) that starts at zero, denoted $\int^{\cdot}_0 x(s) \; \mathrm{d}B^i_s : [0,T] \times \Omega \rightarrow \mathbb{R}^n$. For $\ell \ge2$, we denote by $L^{\ell}_{\mathcal{F}}(\Omega;C([0,T];\mathbb{R}^n))$ the space of $\mathcal{F}$--adapted processes $x : [0,T] \times \Omega \rightarrow \mathbb{R}^n$ that have continuous sample paths and satisfy $\mathbb{E}\left[ \underset{s \in [0,T]}{\sup} \ \| x(s,\omega) \|^{\ell} \right] < \infty$.

\subsection{Stochastic Differential Equations}

From now on, we fix two integers $n, m \in \mathbb{N}$, a maximal time $T \in \mathbb{R}_+$, and a compact, convex subset $U \subseteq \mathbb{R}^m$. Although in this work we consider differential equations steered by deterministic controls (see Section \ref{sec_ProblemFormulation} below), for the sake of generality, we introduce stochastic dynamics which depend on either deterministic or stochastic controls. Specifically, we denote by $\mathcal{U}$  the set of admissible controls and consider either deterministic controls $\mathcal{U} = L^2([0,T];U)$ or stochastic controls $\mathcal{U} = L^2_{\mathcal{F}}([0,T] \times \Omega;U)$. Note that since $U$ is compact, admissible controls are almost everywhere, or a.e. for brevity (and additionally almost surely, or a.s. for brevity) bounded. We are given continuous mappings $b_i : \mathbb{R} \times \mathbb{R}^n \rightarrow \mathbb{R}^n$, $i=0,\dots,m$ and $\sigma_j : \mathbb{R} \times \mathbb{R}^n \rightarrow \mathbb{R}^n$, $j=1,\dots,d$ which are at least $C^2$ with respect to the variable $x$. For a given $u \in \mathcal{U}$, we consider dynamical systems modeled through the following forward stochastic differential equation with uncontrolled diffusion
\begin{align} \label{eq_SDE}
    \displaystyle \mathrm{d}x(t) &= b(t,u(t),x(t)) \; \mathrm{d}t + \sigma(t,x(t)) \; \mathrm{d}B_t , \quad x(0) = x^0 \\
    \displaystyle &\triangleq \left( b_0(t,x(t)) + \sum^m_{i=1} u^i(t) b_i(t,x(t)) \right) \; \mathrm{d}t + \sum^d_{j=1} \sigma_j(t,x(t)) \; \mathrm{d}B^j_t \nonumber
\end{align}
where we assume that the fixed initial condition satisfies $x^0 \in L^{\ell}_{\mathcal{F}_0}(\Omega;\mathbb{R}^n)$, for every $\ell \ge 2$ (for instance, this holds when $x^0$ is a deterministic vector of $\mathbb{R}^n$).

The procedure developed in this work is based on the following linearization of \eqref{eq_SDE}. For $\ell \ge 2$, let $v \in \mathcal{U}$ and $y \in L^{\ell}_{\mathcal{F}}(\Omega;C([0,T];\mathbb{R}^n))$. For a given $u \in \mathcal{U}$, we define the linearization of \eqref{eq_SDE} around $(v,y)$ to be the following well-defined forward stochastic differential equation with uncontrolled diffusion
{
\begin{align} \label{eq_LSDE}
    \displaystyle &\mathrm{d}x(t) = b_{v,y}(t,u(t),x(t)) \; \mathrm{d}t + \sigma_y(t,x(t)) \; \mathrm{d}B_t , \quad x(0) = x^0 \\
    \displaystyle &\triangleq \left( b_0(t,y(t)) + \frac{\partial b_0}{\partial x}(t,y(t)) (x(t) - y(t)) + \sum^m_{i=1} \left( u^i(t) b_i(t,y(t)) + v^i(t) \frac{\partial b_i}{\partial x}(t,y(t)) (x(t) - y(t)) \right) \right) \mathrm{d}t \nonumber \\
    \displaystyle &+ \sum^d_{j=1} \left( \sigma_j(t,y(t)) + \frac{\partial \sigma_j}{\partial x}(t,y(t)) (x(t) - y(t)) \right) \mathrm{d}B^j_t . \nonumber
\end{align}}

We require the solutions to \eqref{eq_SDE} and to \eqref{eq_LSDE} to be bounded in expectation uniformly with respect to $u$, $v$ and $y$. For this, we consider the following (standard) assumption:

\vspace{5pt}

\noindent $(A_1)$ Functions $b_i$, $i=0,\dots,m$, $\sigma_j$, $j=1,\dots,d$, have compact supports in $[0,T] \times \mathbb{R}^n$.

\vspace{5pt}
    
\noindent Under $(A_1)$, for every $\ell \ge 2$ and every $u \in \mathcal{U}$, the stochastic equation \eqref{eq_SDE} has a unique (up to stochastic indistinguishability) solution $x_u \in L^{\ell}_{\mathcal{F}}(\Omega;C([0,T];\mathbb{R}^n))$, whereas for every $u, v \in \mathcal{U}$ and $y \in L^{\ell}_{\mathcal{F}}(\Omega;C([0,T];\mathbb{R}^n))$, the stochastic equation \eqref{eq_LSDE} has a unique (up to stochastic indistinguishability) solution $x_{u,v,y} \in L^{\ell}_{\mathcal{F}}(\Omega;C([0,T];\mathbb{R}^n))$ (see, e.g., \cite{LeGall2016,Carmona2016}), and the following technical result holds, which will be used in the proof of our convergence result, with proof given in the appendix (see Section \ref{sec_app1}).

\begin{lmm}[Uniform boundness and continuity with respect to control inputs] \label{lemma_bound}
Fix $\ell \ge 2$, and for $u \in \mathcal{U}$, let $x_u$ denote the solution to \eqref{eq_SDE}, whereas for $u, v \in \mathcal{U}$ and $y \in L^{\ell}_{\mathcal{F}}(\Omega;C([0,T];\mathbb{R}^n))$, let $x_{u,v,y}$ denote the solution to \eqref{eq_LSDE}. Under $(A_1)$, there exists a constant $C \ge 0$ which does not depend on $u$, $v$, or $y$ such that it holds that
$$
\mathbb{E}\left[ \underset{t \in [0,T]}{\sup} \  \| x_u(t) \|^{\ell} + \underset{t \in [0,T]}{\sup} \  \| x_{u,v,y}(t) \|^{\ell} \right] \le C ,
$$
{\small
$$\displaystyle \mathbb{E}\left[ \underset{t \in [0,T]}{\sup} \ \| x_{u_1}(t) - x_{u_2}(t) \|^{\ell} + \underset{t \in [0,T]}{\sup} \ \| x_{u_1,v,y}(t) - x_{u_2,v,y}(t) \|^{\ell} \right] \le C \mathbb{E}\left[ \left( \int^T_0 \| u_1(s) - u_2(s) \| \; \mathrm{d}s \right)^{\ell} \right], u_1, u_2 \in \mathcal{U} .
$$}
\end{lmm}

\subsection{Stochastic Optimal Control Problem} \label{sec_ProblemFormulation}

Given $q \in \mathbb{N}$, we consider continuous mappings $g : \mathbb{R}^n \rightarrow \mathbb{R}^q$, $G : \mathbb{R}^m \rightarrow \mathbb{R}$, $H : \mathbb{R}^n \rightarrow \mathbb{R}$, and $L : \mathbb{R} \times \mathbb{R}^m \times \mathbb{R}^n \rightarrow \mathbb{R}$ with $L(t,x,u) \triangleq L_0(t,x) + \sum^m_{i=1} u^i L_i(t,x)$. We require $g$, $H$, and $L_i$, $i=0,\dots,m$, to be at least $C^2$ with respect to the variable $x$, and we require $G$, $H$ to be convex. In particular, $G$ is Lipschitz when restricted to the compact and convex set $U$. We focus on finite-horizon, finite-dimensional non-linear stochastic Optimal Control Problems (OCP) with control-affine dynamics and uncontrolled diffusion, of the form
{%
\begin{equation*}
    \begin{cases}
        \displaystyle \underset{u \in \mathcal{U}}{\min} \ \mathbb{E}\left[ \int^T_0 f^0(s,u(s),x(s)) \; \mathrm{d}s \right] \triangleq \mathbb{E}\left[ \int^T_0 \Big( G(u(s)) + H(x(s)) + L(s,u(s),x(s)) \Big) \; \mathrm{d}s \right] \medskip \\
        \mathrm{d}x(t) = b(t,u(t),x(t)) \; \mathrm{d}t + \sigma(t,x(t)) \; \mathrm{d}B_t , \quad x(0) = x^0 , \quad \mathbb{E}\left[ g(x(T)) \right] = 0
    \end{cases}
\end{equation*}}
where we optimize over deterministic controls $u \in \mathcal{U} = L^2([0,T];U)$. We adopt the (fairly mild) assumption:

\vspace{5pt}

\noindent $(A_2)$ Mappings $g$, $H$, and $L_i$, $i=0,\dots,m$, either are affine-in-state or have compact supports in $\mathbb{R}^n$ and in $[0,T] \times \mathbb{R}^n$, respectively.

\vspace{5pt}

Our choice of optimizing over deterministic controls as opposed to stochastic controls is motivated by practical considerations. Specifically, in several applications of interest ranging from aerospace to robotics, it is often advantageous to compute and implement simpler deterministic controls at higher control rates, to be able to quickly react to external disturbances, unmodeled dynamical effects, and changes in the cost function (e.g., moving obstacles that a robot should avoid in real-time). 
Moreover, in cases where a feedback controller is accounted for, a common and efficient approach entails decomposing the stochastic control into a state-dependent feedback term, plus a nominal deterministic control trajectory to be optimized for, which is equivalent to adopting deterministic controls (see, e.g., \cite{Kazuhide2018,LewBonalli2020}). Nevertheless, for the sake of completeness and generality, in Section \ref{sec_extension} we introduce \firstRev{appropriate} conditions under which our method may extend to stochastic controls; we also analyze possible extensions to the case of free-final-time optimal control problems.

Many applications of interest often involve state constraints. In this case, to make sure the procedure developed in this work still applies, every such constraint needs to be considered in expectation and penalized within the cost of OCP (for example by including those contributions in $H$ or $L$ through some penalization function). \firstRev{In future work,} we plan to extend our method to more general settings, whereby for instance state constraints are enforced through more accurate chance constraints (see Section \ref{sec_conclusion} for a thorough discussion).

\section{Stochastic Sequential Convex Programming} \label{sec_SCP}

We propose the following framework to solve OCP, based on the classical SCP methodology. Starting from some initial guesses of control $u_0 \in \mathcal{U}$ and trajectory $x_0 \in L^{\ell}_{\mathcal{F}}(\Omega;C([0,T];\mathbb{R}^n))$, $\ell \ge 2$, we inductively define a sequence of stochastic linear-convex problems whose dynamics and costs stem from successive linearizations of the mappings $b$, $\sigma$, and $L$, and we successively solve those problems while updating user-defined parameters. \firstRev{The convergence of the method generally depends on} a good choice of the initial guess $(u_0,x_0)$ and of updates in each iteration to \textit{trust-region constraints}. These constraints are added to make the successive linearizations of OCP well-posed. Below, we detail this procedure.

\subsection{The Method}

At iteration $k+1 \in \mathbb{N}$, by denoting
\begin{equation} \label{eq_LCost}
    f^0_{v,y}(s,u,x) \triangleq G(u) + H(x) + L(s,u,y) + \frac{\partial L}{\partial x}(s,v,y) (x - y) ,
\end{equation}
we define the following stochastic Linearized Optimal Control Problem (LOCP$^{\Delta}_{k+1}$)
\begin{equation*}
    \begin{cases}
        \displaystyle \underset{u \in \mathcal{U}}{\min} \ \mathbb{E}\left[ \int^T_0 f^0_{k+1}(s,u(s),x(s)) \; \mathrm{d}s \right] \triangleq \mathbb{E}\bigg[ \int^T_0 f^0_{u_k,x_k}(s,u(s),x(s)) \; \mathrm{d}s \bigg] \medskip \\
        \mathrm{d}x(t) = b_{k+1}(t,u(t),x(t)) \; \mathrm{d}t + \sigma_{k+1}(t,x(t)) \; \mathrm{d}B_t , \quad x(0) = x^0 \medskip \\
        \hspace{6ex} \triangleq b_{u_k,x_k}(t,u(t),x(t)) \; \mathrm{d}t + \sigma_{x_k}(t,x(t)) \; \mathrm{d}B_t \medskip \\
        \displaystyle \mathbb{E}\left[ g_{k+1}(x(T)) \right] \triangleq \mathbb{E}\left[ g(x_k(T)) + \frac{\partial g}{\partial x}(x_k(T)) (x(T) - x_k(T)) \right] = 0 \medskip \\
        \displaystyle \int^T_0 \mathbb{E}\left[ \| x(s) - x_k(s) \|^2 \right] \; \mathrm{d}s \le \Delta_{k+1}
    \end{cases}
\end{equation*}
where we optimize over deterministic controls $u \in \mathcal{U} = L^2([0,T];U)$. The tuple $(u_k,x_k) \in \mathcal{U} \times L^{\ell}_{\mathcal{F}}(\Omega;C([0,T];\mathbb{R}^n))$ is defined inductively and denotes a solution to (LOCP)$^{\Delta}_k$.

Each problem LOCP$^{\Delta}_{k+1}$ consists of linearizing OCP around the \textit{solution at the previous iteration} $(u_k,x_k)$, starting from $(u_0,x_0)$. To avoid misguidance due to high linearization error, we must restrict the search for optimal solutions for LOCP$^{\Delta}_{k+1}$ to neighborhoods of $(u_k,x_k)$. This is achieved by the final constraints listed in LOCP$^{\Delta}_{k+1}$, referred to as trust-region constraints, where the constant $\Delta_{k+1} \ge 0$ is the trust-region radius. No such constraint is enforced on controls, since those appear linearly in $b$, $\sigma$, and $L$. \firstRev{To ensure well-posedness of} the sequence (LOCP$^{\Delta}_k$)$_{k \in \mathbb{N}}$, we require LOCP$^{\Delta}_{k+1}$ has a solution\firstRev{, for which we consider the following assumption:}

\vspace{5pt}

\noindent $(A_3)$ For every $k \in \mathbb{N}$, problem LOCP$^{\Delta}_{k+1}$ is feasible.

\begin{prpstn}[Existence of solutions of LOCP$^{\Delta}_{k+1}$]
\firstRev{Under $(A_3)$, LOCP$^{\Delta}_{k+1}$ has a solution for every $k \in \mathbb{N}$.}
\end{prpstn}

\begin{proof}
\firstRev{The proof is based on the argument developed to prove \cite[Theorem 5.2, Chapter 2]{Yong1999}. Specifically, let $k \in \mathbb{N}$ and $\{ (u^{\alpha}_{k+1},x^{\alpha}_{k+1}) \}_{\alpha \in \mathbb{N}} \in \mathcal{U} \times L^{\ell}_{\mathcal{F}}(\Omega;C([0,T];\mathbb{R}^n))$ be a minimizing sequence for LOCP$^{\Delta}_{k+1}$.}

\firstRev{The sequence $\{ u^{\alpha}_{k+1} \}_{\alpha \in \mathbb{N}}$ is uniformly bounded in $L^2([0,T];\mathbb{R}^m)$, and thus we may assume the existence of $\tilde u_{k+1} \in L^2([0,T];\mathbb{R}^m)$ such that $\{ u^{\alpha}_{k+1} \}_{\alpha \in \mathbb{N}}$ converges, up to some subsequence, to $\tilde u_{k+1}$ for the weak topology of $L^2([0,T];\mathbb{R}^m)$. Moreover, the compactness and convexity of $U$ yield $\tilde u_{k+1} \in \mathcal{U}$. Finally, by Mazur's theorem there exists a sequence of convex combinations
$$
\tilde u^{\alpha}_{k+1} \triangleq \sum_{\beta \ge 1} c^{k+1}_{\alpha,\beta} u^{\alpha+\beta}_{k+1} , \quad c^{k+1}_{\alpha,\beta} \ge 0 , \quad \sum_{\beta \ge 1} c^{k+1}_{\alpha,\beta} = 1 ,
$$
such that $\{ \tilde u^{\alpha}_{k+1} \}_{\alpha \in \mathbb{N}}$ converges to $\tilde u_{k+1}$ for the strong topology of $L^2([0,T];\mathbb{R}^m)$ (although, since controls are deterministic and the diffusion is uncontrolled, weak convergence of controls would suffice in our setting).}

\firstRev{Thanks to Lemma \ref{lemma_bound}, the sequence of trajectories $\{ x_{\tilde u^{\alpha}_{k+1},u_k,x_k} \}_{\alpha \in \mathbb{N}} \subseteq L^{\ell}_{\mathcal{F}}(\Omega;C([0,T];\mathbb{R}^n))$ converges to the trajectory $\tilde x_{k+1} \triangleq x_{\tilde u_{k+1},u_k,x_k} \in L^{\ell}_{\mathcal{F}}(\Omega;C([0,T];\mathbb{R}^n))$ for the strong topology of $L^{\ell}_{\mathcal{F}}(\Omega;C([0,T];\mathbb{R}^n))$. At this step, for every $\alpha \in \mathbb{N}$, thanks to the linearity of \eqref{eq_LSDE} we obtain that
$$
x_{\tilde u^{\alpha}_{k+1},u_k,x_k} \equiv \sum_{\beta \ge 1} c^{k+1}_{\alpha,\beta} x_{u^{\alpha+\beta}_{k+1},u_k,x_k} = \sum_{\beta \ge 1} c^{k+1}_{\alpha,\beta} x^{\alpha+\beta}_{k+1} ,
$$
and therefore the linearity of the function $g_{k+1}$ yields
$$
\mathbb{E}\left[ g_{k+1}(x_{\tilde u^{\alpha}_{k+1},u_k,x_k}(T)) \right] = \sum_{\beta \ge 1} c^{k+1}_{\alpha,\beta} \mathbb{E}\left[ g_{k+1}(x^{\alpha+\beta}_{k+1}(T)) \right] = 0 ,
$$
whereas the convexity of the norm yields
$$
\left( \int^T_0 \mathbb{E}\left[ \| x_{\tilde u^{\alpha}_{k+1},u_k,x_k}(s) - x_k(s) \|^2 \right] \; \mathrm{d}s \right)^{\frac{1}{2}} \le \sum_{\beta \ge 1} c^{k+1}_{\alpha,\beta} \left( \int^T_0 \mathbb{E}\left[ \| x^{\alpha+\beta}_{k+1}(s) - x_k(s) \|^2 \right] \; \mathrm{d}s \right)^{\frac{1}{2}} \le \Delta^{\frac{1}{2}}_{k+1} .
$$
In turn, passing to the limit for $\alpha \rightarrow \infty$, we infer that the tuple $(\tilde u_{k+1},\tilde x_{k+1}) \in \mathcal{U} \times L^{\ell}_{\mathcal{F}}(\Omega;C([0,T];\mathbb{R}^n))$ is admissible for LOCP$^{\Delta}_{k+1}$. In addition, the convexity of the function $f^0_{k+1}$ allows us to compute
\begin{align*}
    \mathbb{E}\left[ \int^T_0 f^0_{k+1}(s,\tilde u_{k+1}(s),\tilde x_{k+1}(s)) \; \mathrm{d}s \right] &\le \underset{\alpha \rightarrow \infty}{\lim} \ \sum_{\beta \ge 1} c^{k+1}_{\alpha,\beta} \mathbb{E}\left[ \int^T_0 f^0_{k+1}(s,u^{\alpha + \beta}_{k+1}(s),x^{\alpha+\beta}_{k+1}(s)) \; \mathrm{d}s \right] \\
    &= \underset{u \in \mathcal{U}}{\min} \ \mathbb{E}\left[ \int^T_0 f^0_{k+1}(s,u(s),x(s)) \; \mathrm{d}s \right] ,
\end{align*}
from which we conclude that $(\tilde u_{k+1},\tilde x_{k+1}) \in \mathcal{U} \times L^{\ell}_{\mathcal{F}}(\Omega;C([0,T];\mathbb{R}^n))$ is a solution of LOCP$^{\Delta}_{k+1}$.}
\end{proof}

\firstRev{\begin{rmrk}
    Although 
    one finds empirically that $(A_3)$ is often satisfied in practice (see for instance results in Section \ref{sec_numerical}), 
    it is generally difficult to derive sufficient conditions for this assumption to hold true a priori. In particular, to the best of our knowledge, only stochastic dynamics with constant coefficients and with specific regularity conditions on the stochastic diffusion yield controllability \cite{wang2017}. In deterministic settings, some works obtain feasibility of each convex subproblem, though only locally around the unknown solution of the original formulation by assuming second-order regularity conditions, which however can not be checked a priori \cite{Dinh2010,Diehl2019}. Thus, even SCP-based schemes in simpler deterministic settings must often either explicitly assume the feasibility of each convex subproblem, or must modify the linearized dynamics by infusing additional slack controls to force feasibility \cite{Mao2017}. Motivated by these remarks, we leave the investigation for sufficient conditions for the validity of $(A_3)$ as a direction of future research, which are out of the scope of this work which focuses on the properties of accumulations points of the sequence generated by SCP. Note that $(A_3)$ is satisfied when $g = 0$, i.e., no final constraints are imposed: a problem that remains generally relevant though computationally difficult to solve.
\end{rmrk}}

Under assumptions $(A_1)$--$(A_3)$, the method consists of iteratively solving \firstRev{the aforementioned} linearized problems through the update of the sequence of trust-region radii, producing a sequence of tuples $(u_k,x_k)_{k \in \mathbb{N}}$ such that for each $k \in \mathbb{N}$, $(u_{k+1},x_{k+1})$ solves LOCP$^{\Delta}_{k+1}$. The user may \firstRev{often} steer this procedure to convergence (with respect to appropriate topologies) by adequately selecting an initial guess $(u_0,x_0)$ and an update rule for $(\Delta_k)_{k \in \mathbb{N}}$\footnote{When state constraint penalization is adopted, SCP procedures may also consider update rules for \textit{penalization weights}, and those must be provided together with update rules for the trust-region radii. Further details can be found in \cite{palacios1982,boggs1995}.}, and appropriate choices will be described in Section \ref{sec_numerical}. Assuming that an accumulation point for $(u_k,x_k)_{k \in \mathbb{N}}$ can be found (whose existence will be discussed shortly), our objective consists of proving that this is a candidate locally-optimal solution to the original formulation OCP. Specifically, we show that any accumulation point for $(u_k,x_k)_{k \in \mathbb{N}}$ satisfies the stochastic PMP related to OCP. To develop such analysis, we require the absence of state constraints, and in particular of trust-region constraints. 

\subsection{Stochastic Pontryagin Maximum Principle}

\firstRev{In this section, we recall the statement of the PMP which provides classical first-order necessary conditions for optimality, upon which we will establish our main result. For the sake of clarity, we introduce the PMP related to OCP and the PMP related to each convexified problem LOCP$^{\Delta}_k$ separately}.

\subsubsection{\firstRev{PMP related to OCP}}

\firstRev{For every $p \in \mathbb{R}^n$, $p^0 \in \mathbb{R}$, and $q = (q_1,\dots,q_d) \in \mathbb{R}^{n \times d}$ define the \textit{Hamiltonian} (same notation as in \eqref{eq_SDE})}
$$
\firstRev{H(s,u,x,p,p^0,q) \triangleq p^{\top} b(s,x,u) + p^0 f^0(s,x,u) + \sum^d_{i=1} q^{\top}_i \sigma_i(s,x) .}
$$
\begin{thrm}[Stochastic Pontryagin Maximum Principle for OCP \cite{Yong1999}] \label{Theo_PMP}
Let $(u,x)$ be a locally-optimal solution to OCP. There exist $p\in L^2_{\mathcal{F}}(\Omega;C([0,T];\mathbb{R}^n))$, a tuple $(\mathfrak{p},p^0)$, where $\mathfrak{p} \in \mathbb{R}^d$ and $p^0 \le 0$ are constant, and $q = (q_1,\dots,q_d) \in L^2_{\mathcal{F}}([0,T] \times \Omega;\mathbb{R}^{n \times d})$ such that the following relations are satisfied:
\begin{enumerate}
    \item Non-Triviality Condition: $(\mathfrak{p},p^0) \neq 0$.
    \item Adjoint Equation:
    \firstRev{\begin{align*}
        \mathrm{d}p(t) &= \displaystyle -\frac{\partial H}{\partial x}(t,u(t),x(t),p(t),p^0,q(t)) \; \mathrm{d}t + q(t) \; \mathrm{d}B_t , \quad p(T) = \displaystyle \mathbb{E}\left[ \frac{\partial g}{\partial x}(x(T)) \right]^{\top} \mathfrak{p} \in \mathbb{R}^n .
    \end{align*}}
    \item Maximality Condition:
    $$
    \firstRev{u(t) = \underset{v \in U}{\arg \max} \ \mathbb{E}\Big[ H(t,v,x(t),p(t),p^0,q(t)) \Big] , \ \textnormal{a.e.}}
    $$
\end{enumerate}
\firstRev{The quantity $(u,\mathfrak{p},p^0)$ uniquely determines $x$, $p$, and $q$ and is called extremal for OCP (associated with the tuple $(u,x,p,\mathfrak{p},p^0,q)$, or simply with $(u,x)$). An extremal $(u,\mathfrak{p},p^0)$ is called \textnormal{normal} if $p^0 \neq 0$}.
\end{thrm}

Albeit final conditions are specified instead of initial conditions for $p$, it turns out that processes satisfying backward stochastic differential equations are adapted with respect to the filtration $\mathcal{F}$ (see, e.g., \cite{Carmona2016,Yong1999}), which makes the adjoint equation well-posed. Although conditions for optimality for stochastic optimal control problems are usually developed when considering stochastic controls only, the proof of Theorem \ref{Theo_PMP} \firstRev{(and of Theorem \ref{Theo_PMPk} below)} %
readily follows from classical arguments (e.g., see \cite[Chapter 3.6]{Yong1999}). Nevertheless, to prove our main result, we rely on a proof of Theorem \ref{Theo_PMP} \firstRev{(and of Theorem \ref{Theo_PMPk} below)} which stems from implicit-function-theorem-type results (see Sections \ref{sec_PMP} and \ref{sec_proofConv}). Thus, in Section \ref{sec_PMP}, we provide a new proof of Theorem \ref{Theo_PMP} \firstRev{(more precisely, of Theorem \ref{Theo_PMPk} below)} which follows from the original idea developed by Pontryagin and his group (see, e.g., \cite{pontryagin2018,gamkrelidze2013,agrachev2013}), though the latter result should not be understood as part of our main contribution.

\subsubsection{\firstRev{PMP related to LOCP$^{\Delta}_k$}}
\firstRev{For every $k \ge 1$, $p \in \mathbb{R}^n$,  $p^0 , p^1 \in \mathbb{R}$, and $q = (q_1,\dots,q_d) \in \mathbb{R}^{n \times d}$, define the \textit{Hamiltonian} (same notation as in \eqref{eq_LSDE})
$$
H_k(s,u,x,p,p^0,p^1,q) \triangleq p^{\top} b_k(s,x,u) + p^0 f^0_k(s,x,u) + p^1 \| x - x_{k-1}(s) \|^2 + \sum^d_{i=1} q^{\top}_i (\sigma_k)_i(s,x) .
$$}
\begin{thrm}[Weak Stochastic Pontryagin Maximum Principle for LOCP$^{\Delta}_k$ \cite{Yong1999}] \label{Theo_PMPk}
\firstRev{Let $(u_k,x_k)$ be a locally-optimal solution to LOCP$^{\Delta}_k$. There exists $p_k \in L^2_{\mathcal{F}}(\Omega;C([0,T];\mathbb{R}^n))$, a tuple $(\mathfrak{p}_k,p^0_k,p^1_k)$, where $\mathfrak{p}_k \in \mathbb{R}^q$, $p^0_k \le 0$, and $p^1_k \in \mathbb{R}$ are constant, and $q_k = ((q_k)_1,\dots,(q_k)_d) \in L^2_{\mathcal{F}}([0,T] \times \Omega;\mathbb{R}^{n \times d})$ such that the following relations are satisfied:
\begin{enumerate}
    \item Non-Triviality Condition: $(\mathfrak{p}_k,p^0_k,p^1_k) \neq 0$.
    \item Adjoint Equation:
    \begin{align*}
        \mathrm{d}p_k(t) &= \displaystyle -\frac{\partial H_k}{\partial x}(t,u_k(t),x_k(t),p_k(t),p^0_k,p^1_k,q_k(t)) \; \mathrm{d}t + q_k(t) \; \mathrm{d}B_t , \quad p_k(T) = \displaystyle \mathbb{E}\left[ \frac{\partial g_k}{\partial x}(x_k(T)) \right]^{\top} \mathfrak{p}_k \in \mathbb{R}^n .
    \end{align*}
    \item Maximality Condition:
    $$
    u_k(t) = \underset{v \in U}{\arg \max} \ \mathbb{E}\Big[ H_k(t,v,x_k(t),p_k(t),p^0_k,p^1_k,q_k(t)) \Big] , \ \textnormal{a.e.}
    $$
\end{enumerate}
The quantity $(u_k,\mathfrak{p}_k,p^0_k,p^1_k)$ uniquely determines $x_k$, $p_k$, and $q_k$ and is called extremal for LOCP$^{\Delta}_k$ (associated with $(u_k,x_k,p_k,\mathfrak{p}_k,p^0_k,p^1_k,q_k)$, or simply with $(u_k,x_k)$). An extremal $(u_k,\mathfrak{p}_k,p^0_k,p^1_k)$ is called \textnormal{normal} if $p^0_k \neq 0$}.
\end{thrm}

\firstRev{\begin{rmrk} \label{remark_slack}
    By introducing the new variable
    \begin{equation} \label{eq_changeVar}
        y(t) \triangleq \int^t_0 \| x(s) - x_k(s) \|^2 \; \mathrm{d}s ,
    \end{equation}
    the trust-region constraint in LOCP$^{\Delta}_{k+1}$ can be rewritten as $\mathbb{E}[ y(T) - \Delta_{k+1} ] \le 0$. Thus, by leveraging this transformation, LOCP$^{\Delta}_{k+1}$ may be reformulated as a standard stochastic optimal control problem with final inequality constraints. Nevertheless, although the conditions listed in Theorem \ref{Theo_PMP} are essentially sharp, the statement of Theorem \ref{Theo_PMPk} may be strengthened as follows. If $(u_{k+1},x_{k+1},\mathfrak{p}_{k+1},p^0_{k+1},p^1_{k+1},q_{k+1})$ is an extremal for LOCP$^{\Delta}_{k+1}$, one can additionally prove that $p^1_{k+1} \le 0$ (e.g., see \cite{Frankowska2018}; note that in \cite{Frankowska2018} the multipliers have opposite signs because a different convention is adopted) and
    $$
    p^1_{k+1} \mathbb{E}\left[ \int^T_0 \| x_{k+1}(s) - x_k(s) \|^2 \; \mathrm{d}s - \Delta_{k+1} \right] = 0
    $$
    (the latter is know as \textit{slack condition}), which motivates the choice ``weak stochastic PMP for LOCP$^{\Delta}_{k+1}$'' as name for Theorem \ref{Theo_PMPk}. Nevertheless, since the trust-region constraints do not appear in the original problem, we do not need to leverage these latter additional conditions on $p^1_{k+1}$ to prove our claims, i.e., Theorem \ref{Theo_PMPk} suffices to establish the aforementioned properties of accumulation points for SCP when applied to solve OCP.
\end{rmrk}}

\subsection{Main Results} \label{sec_mainDiscussion}

\firstRev{Our contribution is twofold. First, under very mild assumptions, we prove that any accumulation point of the sequence of iterates generated by SCP is a candidate locally-optimal solution for OCP. Specifically, we prove that any accumulation point of the sequence $(u_k,x_k,\mathfrak{p}_k,p^0_k,p^1_k,q_k)_{k \in \mathbb{N}}$ generated by SCP, where $(u_k,\mathfrak{p}_k,p^0_k,p^1_k)$ is an extremal for LOCP$^{\Delta}_k$ associated with $(u_k,x_k)$ in the sense of Theorem \ref{Theo_PMP}, is an extremal for OCP in the sense of Theorem \ref{Theo_PMPk} (see Theorem \ref{theo_main}). Although the optimization community is aware that establishing the convergence of the sequence generated by SCP is generally difficult (see, e.g., \cite{Nocedal1999,Lu2013}), one can often prove optimality-related properties of its accumulation points, a much natural property which justify the use of SCP (see, e.g., \cite{Lu2013}). Second, by strengthening our original assumptions, we prove existence of at least one accumulation point of the aforementioned sequence $(u_k,x_k,\mathfrak{p}_k,p^0_k,p^1_k,q_k)_{k \in \mathbb{N}}$ generated by SCP (see Theorem \ref{theo_main2}). Below, we organize these claims in two more precise statements.}

\subsubsection{\firstRev{Properties of Accumulation Points for SCP}}

\begin{thrm}[Properties of Accumulation Points for SCP] \label{theo_main}
\firstRev{Assume that $(A_1)$--$(A_3)$ hold and that SCP generates a sequence $(\Delta_k,u_k,x_k)_{k \in \mathbb{N}}$ such that $(\Delta_k)_{k \in \mathbb{N}} \subseteq \mathbb{R}_+ \setminus \{ 0 \}$ converges to zero, and for every $k \ge 1$, the tuple $(u_k,x_k)$ locally solves LOCP$^{\Delta}_k$. For every $k \ge 1$, letting $(u_k,\mathfrak{p}_k,p^0_k,p^1_k)$ be an extremal associated with $(u_k,x_k)$ for LOCP$^{\Delta}_k$ (whose existence is ensured by Theorem \ref{Theo_PMPk}), assume the following Accumulation Condition holds:
\begin{itemize}
    \item[\textnormal{(AC)}] Up to some subsequence, $(u_k,\mathfrak{p}_k,p^0_k,p^1_k)$ converges to some $(u,\mathfrak{p},p^0,p^1) \in L^2([0,T];\mathbb{R}^m) \times \mathbb{R}^{q+2}$ for the weak topology of $L^2([0,T];\mathbb{R}^m) \times \mathbb{R}^{q+2}$.
\end{itemize}
If $(\mathfrak{p},p^0) \neq 0$, then $(u,\mathfrak{p},p^0)$ is an extremal for OCP associated with $(u,x_u)$.}
\end{thrm}

\firstRev{The guarantees offered by Theorem \ref{theo_main} read as follows. Under $(A_1)$--$(A_3)$ and by selecting a \textit{shrinking-to-zero} sequence of trust-region radii, if iteratively solving problems LOCP$^{\Delta}_k$ returns a sequence of strategies whose extremals satisfy (AC) with a non-trivial multiplier $(\mathfrak{p},p^0) \neq 0$, then SCP finds a candidate (local) solution to OCP. Theorem \ref{theo_main} extends classical results on the well-posedness of SCP (see, e.g., \cite{Lu2013}) from deterministic to stochastic settings. In particular, the requirement $(\mathfrak{p},p^0) \neq 0$ in Theorem \ref{theo_main} is natural and has an equivalent in deterministic settings, playing the role of some sort of \textit{qualification condition} (see, e.g., \cite[Theorem 3.4]{Lu2013}). Note that the requirement $(\mathfrak{p},p^0) \neq 0$ can be easily numerically checked (see our discussion after Theorem \ref{theo_main2}).}

\subsubsection{\firstRev{Existence of Accumulation Points for SCP}}

Assumptions $(A_1)$--$(A_3)$ together with some minor requirements are sufficient to establish that any accumulation point for the sequence of iterates $(u_k,x_k)_{k \in \mathbb{N}}$ satisfies the stochastic PMP related to OCP \firstRev{(Theorem \ref{theo_main})}. We can additionally infer the existence of accumulation points\firstRev{, i.e., (AC) in Theorem \ref{theo_main} holds true, if some more structure on the data defining OCP is assumed. Importantly, the validity of (AC) endows stochastic SCP with the additional guarantee that at least one accumulation point (with respect to weak topologies) exists, which in turn provides a ``weak'' guarantee of success for the method via the result of Theorem \ref{theo_main} (``weak'' because convergence is satisfied up to some subsequence). We introduce the following technical condition:}

\vspace{5pt}

\noindent \firstRev{$(A_4)$ The mapping $G : \mathbb{R}^{m+1} \rightarrow \mathbb{R}$ is given by $G(t,u) = u^{\top} \mathbb{U}(t) u$, where each $\mathbb{U}(t) \in \mathbb{R}^{m \times m}$ is symmetric definite positive and the mapping $t \mapsto \mathbb{U}(t)^{-1}$ is continuous.}

\vspace{5pt}

\noindent \firstRev{The use we make of $(A_4)$ is essentially contained in the following result.}

\begin{lmm} \label{lemma_continuity}
\firstRev{Under $(A_4)$, for every $k \in \mathbb{N}$, every normal extremal $(u_{k+1},\mathfrak{p}_{k+1},p^0_{k+1},p^1_{k+1})$ for LOCP$^{\Delta}_{k+1}$, i.e., for which $p^0_{k+1} \neq 0$, is such that the corresponding control $u_{k+1}$ is time-continuous.}
\end{lmm}

\begin{proof}
\firstRev{Fix $k \in \mathbb{N}$, and let $(u_{k+1},\mathfrak{p}_{k+1},p^0_{k+1},p^1_{k+1})$ be a normal extremal for LOCP$^{\Delta}_{k+1}$. Due to $p^0_{k+1} \neq 0$, the convexity of $U \subseteq \mathbb{R}^m$ and the maximality condition in Theorem \ref{Theo_PMPk} yield
\begin{equation*}
    u_{k+1}(s) = \left\{
    \begin{aligned}
        \frac{1}{2 p^0_{k+1}} \mathbb{U}^{-1}(s) \gamma_k(s) \quad & \textnormal{if} \quad \frac{1}{2 p^0_{k+1}} \mathbb{U}^{-1}(s) \gamma_k(s) \in U , \\
        \textnormal{Proj}_U\left( \frac{1}{2 p^0_{k+1}} \mathbb{U}^{-1}(s) \gamma_k(s) \right) \quad & \textnormal{if} \quad \frac{1}{2 p^0_{k+1}} \mathbb{U}^{-1}(s) \gamma_k(s) \notin U ,
    \end{aligned}
    \right.
\end{equation*}
where we denote $\gamma_k(s) \triangleq \big( (\gamma_k(s))_1 , \dots , (\gamma_k(s))_m \big)$ with $(\gamma_k(s))_i \triangleq \mathbb{E}\big[ p^{\top}_{k+1}(s) b_i(s,x_k(s)) + p^0_{k+1} L_i(s,x_k(s)) \big]$, whereas $\textnormal{Proj}_U : \mathbb{R}^m \rightarrow \mathbb{R}^m$ denotes the projection over the convex set $U$. The claim readily follows once we prove the mappings $t \in [0,T]  \mapsto \gamma_k(s)$ are continuous, given that $t \in [0,T] \mapsto \mathbb{U}^{-1}(s)$ and $v \mapsto \textnormal{Proj}_U(v)$ are continuous. We only prove the continuity of $t \in [0,T] \mapsto \mathbb{E}[ p_k(t) ]$, given that the continuity of $t \in [0,T]  \mapsto \gamma_k(s)$ can be proved by leveraging similar arguments, and $(A_1)$ and $(A_2)$. For this, thanks to $(A_1)$, $(A_2)$, Lemma \ref{lemma_bound}, Theorem \ref{Theo_PMPk}, and H\"older and Burkholder–Davis–Gundy inequalities, for every $0 \le s < t \le T$, we obtain that
\begin{align*}
    &\mathbb{E}\big[ \| p_k(t) - p_k(s) \| \big] \le \\
    &\le 2 \mathbb{E}\bigg[ \bigg\| \int^t_s \bigg( p^{\top}_k(r) \frac{\partial f_k}{\partial x}(r,x_k(r),u_k(r)) + p^0_k \frac{\partial f^0_k}{\partial x}(r,x_k(r),u_k(r)) + 2 p^1_k ( x_k(r) - x_{k-1}(r) ) \\
    &\hspace{45ex} + \sum^d_{i=1} (q_k)^{\top}_i(r) \frac{\partial (\sigma_k)_i}{\partial x}(r,x_k(r)) \bigg) \; \mathrm{d}r \bigg\| \bigg]  + 2 \mathbb{E}\left[ \left\| \int^t_s q_k(r) \; \mathrm{d}B_r \right\| \right] \\
    &\le C \left( p^0_k + p^1_k + \mathbb{E}\left[ \underset{t \in [0,T]}{\sup} \ \| x_k(t) \|^2 \right]^{\frac{1}{2}} + \mathbb{E}\left[ \underset{t \in [0,T]}{\sup} \ \| p_k(t) \|^2 \right]^{\frac{1}{2}} + \mathbb{E}\left[ \int^T_0 \ \| q_k(t) \|^2 \; \mathrm{d}t \right]^{\frac{1}{2}} \right) ( t - s )^{\frac{1}{2}} ,
\end{align*}
for some appropriate constant $C \ge 0$, and the conclusion follows.}
\end{proof}

\firstRev{\begin{rmrk} \label{remark_data}
    Through Lemma \ref{lemma_continuity}, $(A_4)$ becomes crucial to ensure the validity of (AC) in Theorem \ref{theo_main}. In particular, the works \cite{Haberkorn2011,Bonalli2017,Bonalli2018,bonalli2019}, which analyze continuity properties of extremals with respect to appropriate deformations of some deterministic optimal control problems and which inspired our work, show that the time-continuity of optimal controls $u_k$ for each LOCP$^{\Delta}_k$ represents a requirement which is not easy to relax (in particular, see the counterexample in \cite[Section 2.3]{bonalli2019}), especially in the presence of trust-region constraints. Importantly, motivated by regularity results in deterministic optimal control settings (see, e.g., \cite[Theorem 3.2]{ShvartsmanVinter2006}), we reckon that more generic mappings $G$ might yield Lemma \ref{lemma_continuity}, especially when optimizing over deterministic controls, although we leave the investigation of more general sufficient conditions for the time-continuity of optimal controls $u_k$ for each LOCP$^{\Delta}_k$ as a future research direction, in that it is out of the scope of this work which again focuses on the properties of accumulations points of the sequence generated by SCP.
\end{rmrk}}

\begin{thrm}[Existence of Accumulation Points for SCP] \label{theo_main2}
\firstRev{Assume that $(A_1)$--$(A_4)$ hold and that SCP generates a sequence $(\Delta_k,u_k,x_k)_{k \in \mathbb{N}}$ such that $(\Delta_k)_{k \in \mathbb{N}} \subseteq \mathbb{R}_+ \setminus \{ 0 \}$ converges to zero, and for every $k \ge 1$, the tuple $(u_k,x_k)$ locally solves LOCP$^{\Delta}_k$. For every $k \ge 1$, let $(u_k,\mathfrak{p}_k,p^0_k,p^1_k)$ be an extremal associated with $(u_k,x_k)$ for LOCP$^{\Delta}_k$ (whose existence is ensured by Theorem \ref{Theo_PMPk}). If $p^0_k \neq 0$ for every $k \ge 1$, then (AC) in Theorem \ref{theo_main} holds true.}

\firstRev{In addition, if $(u,\mathfrak{p},p^0)$ denotes the extremal for OCP associated with some $(u,x,p,\mathfrak{p},p^0,q)$, which is provided by Theorem \ref{theo_main}, the following convergence holds, up to some subsequence, for every $\ell \ge 2$ when $k \rightarrow \infty$:
\begin{equation} \label{eq_conv}
    \| (\mathfrak{p}_k,p^0_k) - (\mathfrak{p},p^0) \| + \mathbb{E}\left[ \underset{s \in [0,T]}{\sup} \ \left\| x_k(s) - x(s) \right\|^{\ell} + \underset{s \in [0,T]}{\sup} \ \left\| p_k(s) - p(s) \right\|^2 + \int^T_0 \left\| q_k(s) - q(s) \right\|^2 \; \mathrm{d}s \right] \rightarrow 0 .
\end{equation}}
\end{thrm}

\firstRev{The main takeaway of Theorem \ref{theo_main2} is that an accumulation point for stochastic SCP always exists, which in turn is a candidate (local) solution to OCP due to Theorem \ref{theo_main}, as soon as appropriate qualification-type conditions are satisfied. In particular, similarly to the condition $(\mathfrak{p},p^0) \neq 0$ in Theorem \ref{theo_main}, the requirement $p^0_k \neq 0$ in Theorem \ref{theo_main2} plays the role of an additional qualification condition, but at the level of each subproblem LOCP$^{\Delta}_k$. Note that the latter is a \textit{generic} property in deterministic optimal control (see, e.g., \cite{ChitourJeanEtAl2008}). Finally, the requirement $(\mathfrak{p},p^0) \neq 0$ in Theorem \ref{theo_main} from which (local) optimality of the strategy found by SCP stems can be numerically checked thanks to \eqref{eq_conv}, as soon as the multipliers $(\mathfrak{p}_k,p^0_k)$ are accessible through SCP iterations.}

\subsubsection{\firstRev{Insights to Speed Up the Convergence of SCP}}

In Theorem \ref{theo_main2}, there are also insightful statements concerning the convergence of Pontryagin extremals. Let us outline how those statements may be leveraged to speed up convergence.

For this, adopt the notation of Theorem \ref{Theo_PMP} and assume that we are in the situation where applying the maximality condition of the PMP to problem OCP leads to smooth candidate optimal controls, as functions of the variables $x$, $p^0$, and $p$ (be aware that this might not be straightforward to obtain). We are then in a position to define two-point boundary value problems to solve OCP, also known as \textit{shooting methods}, for which the decision variables become $p^0$, $p(T)$, and $q$. In particular, the core of the method consists of iteratively choosing $(p^0,p(T),q)$ and making the adjoint equation evolve until some given final condition is met (see, e.g., \cite{bryson1975,betts1998} for a more detailed explanation of shooting methods). In the context of deterministic optimal control, when convergence is achieved, shooting methods terminate quite fast (at least quadratically). However, here the bottlenecks are: 1) to deal with the presence of the variable $q$ and 2) to find a good guess for the initial value of $(p^0,p(T),q)$ to make the whole procedure converge. In the setting of Theorem \ref{theo_main2}, a valid option to design well-posed shooting methods is as follows. With the notation and assumptions of Theorem \ref{theo_main2}, up to some subsequence it holds that $(p^0_k,p_k(T),q_k) \rightarrow (p^0,p(T),q)$ (with respect to appropriate topologies) as $k \rightarrow \infty$. Therefore, assuming we have access to Pontryagin extremals along iterations and given some large enough iteration $k$, we can fix $q = q_k$ and initialize with $(p^0_k,p_k(T))$ a shooting method for OCP that operates on the finite-dimensional variable $(p^0,p(T)) \in \mathbb{R}^{q+1}$. If successful, this strategy would speed up the convergence of the entire numerical scheme\firstRev{, though we leave its investigation as a future research direction.}

\section{Proof of the Main Results} \label{sec_proof}

\firstRev{Since the statement of Theorem \ref{theo_main} is in particular contained in Theorem \ref{theo_main2}, it is sufficient to prove Theorem \ref{theo_main2} only. We split this proof into three main steps.} First, we retrace the proof of the stochastic PMP to introduce necessary notation and expressions. In addition, we leverage this step to provide novel insight on how to prove the stochastic PMP by following the lines of the original work of Pontryagin and his group (see, e.g., \cite{pontryagin2018,gamkrelidze2013,agrachev2013}), a proof that we could not find in the stochastic literature. Second, we show the convergence of trajectories and controls, together with the convergence of variational inequalities (see Section \ref{sec_variational} for a definition). The latter represents the cornerstone of the proof and paves the way for the final step, which consists of proving the convergence of the Pontryagin extremals. For the sake of clarity and brevity and without loss of generality, we carry out the proof in the case of scalar Brownian motion, i.e., we assume $d = 1$. Moreover, for any $x \in \mathbb{R}^{n}$ with $n \in \mathbb{N}$, we adopt the notation $\tilde x \triangleq (x,x^{n+1}) \in \mathbb{R}^{n+1}$.

\subsection{Main Steps of the Proof of the Stochastic Maximum Principle} \label{sec_PMP}

\firstRev{For the sake of clarity and brevity, we retrace the proof of Theorem \ref{Theo_PMP} only. The proof of Theorem \ref{Theo_PMPk} follows from a straightforward modification of the steps we provide below, by introducing the additional final constraint $\mathbb{E}[ y(T) - \Delta_{k+1} ] \le 0$ via \eqref{eq_changeVar}. In particular, we highlight those modifications below and in Section \ref{sec_convVariational}.}

\subsubsection{Linear Stochastic Differential Equations}

Define the stochastic matrices $A(t) \triangleq \frac{\partial (b,f^0)}{\partial x}(t,u(t),x(t))$ and $D(t) \triangleq \frac{\partial (\sigma,0)}{\partial x}(t,u(t),x(t))$. For any time $r \in [0,T]$ and any bounded initial condition $\tilde \xi_r \in L^2_{\mathcal{F}_r}(\Omega;\mathbb{R}^{n+1})$, the following problem
\begin{equation} \label{eq_rightSDE}
    \begin{cases}
        \mathrm{d}z(t) = A(t) z(t) \; \mathrm{d}t + D(t) z(t) \; \mathrm{d}B_t \medskip \\
        z(s) = 0 , \ s\in [0,r) , \quad z(r) = \tilde \xi_r
    \end{cases}
\end{equation}
is well-posed \cite{Yong1999}. Its unique solution is the $\mathcal{F}$--adapted with right-continuous sample paths process $z : [0,T] \times \Omega \rightarrow \mathbb{R}^{n+1} : (t,\omega) \mapsto \mathbbm{1}_{[r,T]}(t) \phi(t,\omega) \psi(r,\omega) \xi_r(\omega)$, where the matrix-valued $\mathcal{F}$--adapted with continuous sample paths processes $\phi$ and $\psi$ satisfy
\begin{equation} \label{eq_phiPsi}
    \begin{cases}
        \mathrm{d}\phi(t) = A(t) \phi(t) \mathrm{d}t + D(t) \phi(t) \mathrm{d}B_t \medskip \\
        \phi(0) = I ,
    \end{cases} \qquad \begin{cases}
        \mathrm{d}\psi(t) = -\psi(t) \left( A(t) - D(t)^2 \right) \mathrm{d}t - \psi(t) D(t) \mathrm{d}B_t \medskip \\
        \psi(0) = I ,
    \end{cases}
\end{equation}
respectively. In particular, a straightforward application of the It\^o formula shows that $\phi(t) \psi(t) = \psi(t) \phi(t) = I$, and therefore $\psi(t) = \phi(t)^{-1}$, for every $t \in [0,T]$.

\subsubsection{Needle-like Variations and End-point Mapping} \label{sec_needleLike}

One way to prove the PMP comes from the analysis of specific variations called \textit{needle-like variations} on a mapping called the \textit{end-point mapping}. Those concepts are introduced below in the context of optimization over deterministic controls (see Section \ref{sec_extension} for the generalization of this argument to stochastic controls).

Given an integer $j \in \mathbb{N}$, fix $j$ times $0 < t_1 < \dots < t_j < T$ which are Lebesgue points for $u$, and fix $j$ random variables $u_1,\dots,u_j$ such that $u_i \in U$. For fixed scalars $0 \le \eta_i < t_{i+1} - t_i$, $i=1,\dots,j-1$, and $0 \le \eta_j < T - t_j$, the needle-like variation $\pi = \{ t_i , \eta_i , u_i \}_{i = 1,\dots,j}$ of the control $u$ is defined to be the admissible control $u_{\pi}(t) = u_i$ if $t \in [t_i,t_i+\eta_i]$ and $u_{\pi}(t) = u(t)$ otherwise. Denote by $\tilde x_v$ the solution related to an admissible control $v$ of the augmented system 
\begin{equation} \label{eq_augSys}
    \begin{cases}
        \mathrm{d}x(t) = b(t,v(t),x(t)) \; \mathrm{d}t + \sigma(t,x(t)) \; \mathrm{d}B_t , \quad x(0) = x^0 \medskip \\
        \mathrm{d}x^{n+1}(t) = f^0(t,v(t),x(t)) \; \mathrm{d}t , \hspace{14ex} x^{n+1}(0) = 0
    \end{cases}
\end{equation}
and define the mapping $\tilde g : \mathbb{R}^{n+1} \rightarrow \mathbb{R}^{q+1} : \tilde x \mapsto (g(x),x^{n+1})$. For every fixed time $t \in (t_j,T]$, by denoting $\delta_t \triangleq \min \{ t_{i+1} - t_i , t - t_j : i = 1,\dots,j-1 \} > 0$, the end-point mapping at time $t$ is defined to be the function
\begin{align} \label{eq_endPoint}
    F^j_t : \ &\mathcal{C}^j_t \triangleq B^j_{\delta_t}(0) \cap \mathbb{R}^j_+ \rightarrow \mathbb{R}^{q+1} \\
    &(\eta_1,\dots,\eta_j) \mapsto \mathbb{E}\left[ \tilde g(\tilde x_{u_{\pi}}(t)) - \tilde g(\tilde x_u(t)) \right] \nonumber
\end{align}
where $B^j_{\rho}$ is the open ball in $\mathbb{R}^j$ of radius $\rho > 0$. Due to Lemma \ref{lemma_bound}, it is not difficult to see that $F^j_t$ is Lipschitz (see also the argument developed to prove Lemma \ref{lemma_needle} below). In addition, this mapping may be Gateaux differentiated at zero along admissible directions of the cone $\mathcal{C}^j_t$. For this, denote $\tilde b = (b^{\top},f^0)^{\top}$, $\tilde \sigma = (\sigma^{\top},0)^{\top}$ and let $z_{t_i,u_i}$ be the unique solution to \eqref{eq_rightSDE} with $\xi_{t_i} = \tilde b(t_i,u_i,x(t_i)) - \tilde b(t_i,u(t_i),x(t_i))$.

\firstRev{\begin{rmrk}
    For the proof of Theorem \ref{Theo_PMPk}, the only change compared to the proof of Theorem \ref{Theo_PMP} which is required up to this point consists of replacing the function $\tilde g : \mathbb{R}^{n+1} \rightarrow \mathbb{R}^{q+1}$ which defines the end-point mapping \eqref{eq_endPoint} by the function
    $$
    (\tilde g_k,g^{n+2}_k) : \mathbb{R}^{n+2} \rightarrow \mathbb{R}^{q+2} : (\tilde x,y) \mapsto (g_k(x),x^{n+1},y) .
    $$
    Note that this change is consistent since we will effectively make use of the mapping $F^j_t$ at the time $t = T$ only.
\end{rmrk}}

\begin{lmm}[Stochastic needle-like variation formula] \label{lemma_needle}
Let $(\eta_1,\dots,\eta_j) \in \mathcal{C}^j_t$. For any $t > t_j$, it holds that
\begin{align*}
    \bigg\| \mathbb{E}\bigg[ \tilde g(\tilde x_{u_{\pi}}(t)) - \tilde g(\tilde x_u(t)) - \sum^j_{i=1} \eta_i \frac{\partial \tilde g}{\partial \tilde x}(\tilde x_u(t)) z_{t_i,u_i}(t) \bigg] \bigg\| = o\left( \sum^j_{i=1} \eta_i \right) .
\end{align*}
\end{lmm}
The proof of this result is technical (it requires an intense use of stochastic inequalities) but not difficult. We provide an extensive proof of Lemma \ref{lemma_needle} in the appendix in a more general context (see also Section \ref{sec_extension}).

\subsubsection{Variational Inequalities} \label{sec_variational}

The main step in the proof of the PMP goes by contradiction, leveraging Lemma \ref{lemma_needle}. To this end, for every $j \in \mathbb{N}$, define the linear mapping
$$
dF^j_T : \mathbb{R}^j_+ \rightarrow \mathbb{R}^{q+1} : (\eta_1,\dots,\eta_j) \mapsto \sum_{i=1}^j \eta_i \mathbb{E}\left[ \frac{\partial \tilde g}{\partial \tilde x}(\tilde x_u(T)) z_{t_i,u_i}(T) \right] ,
$$
which due to Lemma \ref{lemma_needle}, satisfies
$$
\underset{\alpha > 0 , \alpha \rightarrow 0}{\lim} \ \frac{F^j_T\left( \alpha \eta \right)}{\alpha} = dF^j_T(\eta) ,
$$
for every $\eta \in \mathbb{R}^j_+$. Finally, consider the closed, convex cone of $\mathbb{R}^{q+1}$ given by
\begin{align*}
    K \triangleq \textnormal{Cl} \bigg( \textnormal{Cone} \ \bigg\{ \mathbb{E}\left[ \frac{\partial \tilde g}{\partial \tilde x}(\tilde x_u(T)) z_{t_i,u_i}(T) \right] : \ \textnormal{for} \ u_i \in U \ \textnormal{and} \ t_i \in (0,T) \ \textnormal{is Lebesgue for} \ u \bigg\} \bigg) .
\end{align*}
If $K = \mathbb{R}^{q+1}$, it would hold $dF^j_T(\mathbb{R}^j_+) = K = \mathbb{R}^{q+1}$, and by applying \cite[Lemma 12.1]{agrachev2013}, one would find that the origin is an interior point of $F^j_T(\mathcal{C}^j_T)$. This would imply that $(u,x)$ cannot be optimal for OCP, which gives a contradiction.

The argument above (together with an application of the separation plane theorem) provides the existence of a non-zero vector denoted $\tilde{\mathfrak{p}} = (\mathfrak{p}^{\top},\mathfrak{p}^0) \in \mathbb{R}^{q+1}$ such that the following variational inequality holds
{
\begin{equation} \label{eq_variation}
\tilde{\mathfrak{p}}^{\top} \mathbb{E}\left[ \frac{\partial \tilde g}{\partial \tilde x}(\tilde x_u(T)) z_{r,v}(T) \right] \le 0 , \ r \in [0,T] \ \textnormal{is Lebesgue for} \ u , \ v \in U .
\end{equation}}

\firstRev{\begin{rmrk}
    For the proof of Theorem \ref{Theo_PMPk}, the variational inequalities \eqref{eq_variation} are replaced by the following upgraded variational inequalities which hold for $(u_{k+1},x_{k+1})$, solution to LOCP$^{\Delta}_{k+1}$:
    \begin{align} \label{eq_variationUpgraded}
        \left(\begin{array}{c} \tilde{\mathfrak{p}}_{k+1} \\ p^1_{k+1} \end{array}\right)^{\top} \mathbb{E}\bigg[ \bigg(\begin{array}{cc} \displaystyle \frac{\partial \tilde g}{\partial \tilde x}(\tilde x_k(T)) & 0 \\ 0 & 1 \end{array}\bigg) &\bigg(\begin{array}{c} z^k_{r,v}(T) \\ (z^k_{r,v})^{n+2}(T) \end{array}\bigg) \bigg] \le 0 , \\
        &r \in [0,T] \ \textnormal{is Lebesgue for} \ u_k , \ v \in U , \nonumber
    \end{align}
    where each $((z^k_{r,v})^{\top},(z^k_{r,v})^{n+2})^{\top}$ solves an equation similar to \eqref{eq_rightSDE} (we report this new equation in Section \ref{sec_convVariational}).
\end{rmrk}}

\subsubsection{Conclusion of the Proof of the Stochastic Maximum Principle} \label{sec_daulDerivation}

The conditions of the PMP are derived by working out the variational inequality \eqref{eq_variation} and finding expressions of some appropriate conditional expectations. The main details are developed below in the context of optimization over deterministic controls (see Section \ref{sec_extension} for the generalization to stochastic controls).

First, by appropriately developing solutions to \eqref{eq_rightSDE}, \eqref{eq_variation} can be rewritten as
$$
\mathbb{E}\left[ \left( \tilde{\mathfrak{p}}^{\top} \frac{\partial \tilde g}{\partial \tilde x}(\tilde x_u(T)) \phi(T) \psi(r) \right)^{\top} \left( \left(\begin{array}{c}
    b \\
    f^0
\end{array}\right)(r,v,x(r)) - \left(\begin{array}{c}
    b \\
    f^0
\end{array}\right)(r,u(r),x(r)) \right) \right] \le 0
$$
for every $r \in [0,T]$ Lebesgue point for $u$ and every $v \in U$. Second, again from the structure of \eqref{eq_rightSDE}, it can be readily checked that by denoting
{
\begin{equation} \label{eq_adjointVector}
    p(t) \triangleq \left( \tilde{\mathfrak{p}}^{\top} \mathbb{E}\left[ \frac{\partial \tilde g}{\partial \tilde x}(\tilde x_u(T)) \phi(T) \Big| \mathcal{F}_t \right] \psi(t) \right)_{1,\dots,n} , \quad p^0 \triangleq \left( \tilde{\mathfrak{p}}^{\top} \frac{\partial \tilde g}{\partial \tilde x}(\tilde x_u(T)) \phi(T) \psi(t) \right)_{n+1} ,
\end{equation}}
the quantity $p^0$ is constant in $[0,T]$ (in addition, its negativity can be shown through a standard reformulation of problem OCP, as done in \cite[Section 12.4]{agrachev2013}). Notice that the stochastic process $p : [0,T] \times \Omega \rightarrow \mathbb{R}^n$ is by definition $\mathcal{F}$--adapted. The quantities so far introduced allow one to reformulate the inequality above as
\begin{align*}
    \mathbb{E}\bigg[ p(t)^{\top} \Big( b(t,u(t),x(t)) - b(t,v,x(t)) \Big) + p^0 \Big( f^0(t,u(t),x(t)) - f^0(t,v,x(t)) \Big) \bigg] \ge 0
\end{align*}
for every $t \in [0,T]$ Lebesgue point for $u$ and $v \in U$, from which we infer the maximality condition of the PMP.

It remains to show the existence of the process $q \in L_{\mathcal{F}}(\Omega;L^2([0,T];\mathbb{R}^n))$, the continuity of the sample paths of the process $p$, and the validity of the adjoint equation. For this, remark that, due to Jensen inequality and Lemma \ref{lemma_bound}, the martingale $\left( \mathbb{E}\left[ \frac{\partial \tilde g}{\partial \tilde x}(\tilde x_u(T)) \phi(T) \Big| \mathcal{F}_t \right] \right)_{t \in [0,T]}$ is bounded in $L^2$. Hence, the martingale representation theorem provides the existence of a process $\mu \in L^2_{\mathcal{F}}([0,T] \times \Omega;\mathbb{R}^{(q+1) \times (n+1)})$ such that $\mathbb{E}\left[ \frac{\partial \tilde g}{\partial \tilde x}(\tilde x_u(T)) \phi(T) \Big| \mathcal{F}_t \right] = \mathbb{E}\left[ \frac{\partial \tilde g}{\partial \tilde x}(\tilde x_u(T)) \phi(T) \Big| \mathcal{F}_0 \right] + \int^t_0 \mu(s) \; \mathrm{d}B_s \triangleq N + \chi(t)$, where $N \in \mathbb{R}^{(q+1) \times (n+1)}$ is a constant matrix. The definition in \eqref{eq_adjointVector} immediately gives that the sample paths of the process $p$ are continuous. Next, an application of It\^o formula (component-wise) readily shows that the product $\chi \psi$ satisfies, for $t \in [0,T]$,
\begin{align*}
    (\chi \psi)(t) &= \left( \int^t_0 \mu(s) \; \mathrm{d}B_s \right) \psi(t) = \int^t_0 \mu(s) \psi(s) \; \mathrm{d}B_s - \int^t_0 \mu(s) \psi(s) D(s) \; \mathrm{d}s \\
    & \quad - \int^t_0 \chi(s) \psi(s) \left( A(s) - D(s)^2 \right) \; \mathrm{d}s - \int^t_0 \chi(s) \psi(s) D(s) \; \mathrm{d}B_s .
\end{align*}
Denoting $q(t) \triangleq \left( \tilde{\mathfrak{p}}^{\top} \Big( \mu(t) \psi(t) - \left( N \psi(t) + \chi(t) \psi(t) \right) D(t) \Big) \right)_{1,\dots,n}$, the computations above readily give the adjoint equation of the PMP. Those computations may also be leveraged to show $p \in L^2_{\mathcal{F}}(\Omega;C([0,T];\mathbb{R}^n))$ and $q \in L^2_{\mathcal{F}}([0,T] \times \Omega;\mathbb{R}^n)$ (see, e.g., \cite[Section 7.2]{Yong1999}).

\subsection{Proof of the Convergence Result} \label{sec_proofConv}

Here we enter the core of the proof of Theorem \ref{theo_main}. The convergence of trajectories and controls is addressed first. We devote the last two sections to the convergence of variational inequalities and Pontryagin extremals. For the sake of clarity and brevity, we only consider free-final-time problems. From now on, we implicitly assume $(A_1)$--$(A_3)$.

\subsubsection{Convergence of Controls and Trajectories}

Due to $(A_1)$--$(A_3)$, there exists a sequence of tuples $(u_k,x_k)_{k \in \mathbb{N}}$ such that for every $k \in \mathbb{N}$, the tuple $(u_{k+1},x_{k+1})$ solves LOCP$^{\Delta}_{k+1}$. \firstRev{In what follows, we implicitly adopt the reformulation of each problem convexified problem LOCP$^{\Delta}_{k+1}$, which consist of adding $\mathbb{E}[ y(T) - \Delta_{k+1} ] \le 0$ to the final constraints through the new variable $y$ in \eqref{eq_changeVar}}. If $u \in \mathcal{U}$ denotes an admissible control for OCP that fulfills the conditions of Theorem \ref{theo_main}, we denote by $\tilde x : [0,T] \times \Omega \rightarrow \mathbb{R}^{n+1}$ the $\mathcal{F}$--adapted with continuous sample paths process solution to the augmented system \eqref{eq_augSys} related to OCP with control $u$. The following holds.
\begin{lmm}[Convergence of trajectories] \label{lemma_convTraj}
Assume that the sequence $(\Delta_k)_{k \in \mathbb{N}} \subseteq \mathbb{R}_+$ converges to zero. If the sequence $(u_k)_{k \in \mathbb{N}}$ converges to $u$ for the weak topology of $L^2$, then $\mathbb{E}\left[ \underset{t \in [0,T]}{\sup} \ \| x_k(t) - x(t) \|^{\ell} \right] \rightarrow 0$, as $k \rightarrow \infty$, for every $\ell \ge 2$.
\end{lmm}

\begin{proof}
For every $t \in [0,T]$, we have (below, $C \ge 0$ represents some often overloaded appropriate constant)
{
\begin{align} \label{eq_ineqProofTraj1}
    &\| x_{k+1}(t) - x(t) \|^{\ell} \le C \left\| \int^t_0 \Big( b_0(s,x_k(s)) - b_0(s,x(s)) \Big) \; \mathrm{d}s \right\|^{\ell} \\
    & + C \sum^m_{i=1} \left\| \int^t_0 \Big( u^i_{k+1}(s) b_i(s,x_k(s)) - u^i(s) b_i(s,x(s)) \Big) \; \mathrm{d}s \right\|^{\ell} \nonumber \\
    & + C \left\| \int^t_0 \left( \frac{\partial b_0}{\partial x}(s,x_k(s)) + \sum^m_{i=1} u^i_k(s) \frac{\partial b_i}{\partial x}(s,x_k(s)) \right) ( x_{k+1}(s) - x_k(s) ) \; \mathrm{d}s \right\|^{\ell} \nonumber \\
    & + C \left\| \int^t_0 \left( \sigma(s,x_k(s)) - \sigma(s,x(s)) + \frac{\partial \sigma}{\partial x}(s,x_k(s)) ( x_{k+1}(s) - x_k(s) ) \right) \; \mathrm{d}B_s \right\|^{\ell} \nonumber
\end{align}}
Now, we take expectations. For the last term, we compute
{\small
\begin{align*}
    &\mathbb{E}\left[ \left\| \int^t_0 \left( \sigma(s,x_k(s)) - \sigma(s,x(s)) + \frac{\partial \sigma}{\partial x}(s,x_k(s)) ( x_{k+1}(s) - x_k(s) ) \right) \; \mathrm{d}B_s \right\|^{\ell} \right] \le \\
    & \le C \left( \int^t_0 \mathbb{E}\left[ \underset{r \in [0,s]}{\sup} \ \left\| x_{k+1}(r) - x(r) \right\|^{\ell} \right] \; \mathrm{d}s + \mathbb{E}\left[ \int^T_0 \left\| x_{k+1}(s) - x_k(s) \right\|^{1+(\ell-1)} \; \mathrm{d}s \right] \right) \\
    & \le C \int^t_0 \mathbb{E}\left[ \underset{r \in [0,s]}{\sup} \ \left\| x_{k+1}(r) - x(r) \right\|^{\ell} \right] \; \mathrm{d}s + C \left( \int^T_0 \mathbb{E}\left[ \left\| x_{k+1}(s) - x_k(s) \right\|^2 \right] \mathrm{d}s \right)^{\frac{1}{2}} \left( \underset{k \in \mathbb{N}}{\sup} \ \mathbb{E}\left[ \underset{s \in [0,T]}{\sup} \ \| x_k(s) \|^{2(\ell-1)} \right] \right)^{\frac{1}{2}} \\
    & \le C \left( \int^t_0 \mathbb{E}\left[ \underset{r \in [0,s]}{\sup} \ \left\| x_{k+1}(r) - x(r) \right\|^{\ell} \right] \; \mathrm{d}s + \Delta_{k+1} \right)
\end{align*}}
due to H\"older and Burkholder--Davis--Gundy inequalities, and the last inequality comes from Lemma \ref{lemma_bound}. Similar computations can be carried out for the first and third terms of \eqref{eq_ineqProofTraj1}.

To handle the second term of \eqref{eq_ineqProofTraj1}, we proceed as follows. We see that since $\mathbb{E}\left[ \int^T_0 \| b_i(s,x(s)) \|^2 \; \mathrm{d}t \right] < \infty$ implies $\int^T_0 \| b_i(s,x(s)) \|^2 \; \mathrm{d}t < \infty$, $\mathcal{G}$--almost surely, for every fixed $i = 1,\dots,m$ and $t \in [0,T]$, the convergence of the sequence of controls for the weak topology of $L^2$  entails that $\left\| \int^t_0 \left( u^i_k(s) - u^i(s) \right) b_i(s,x(s)) \; \mathrm{d}s \right\| \rightarrow 0$, $\mathcal{G}$--almost surely as $k \rightarrow \infty$. In addition, $(A_1)$ gives $\left\| \int^b_a \left( u^i_k(s) - u^i(s) \right) b_i(s,x(s)) \; \mathrm{d}s \right\| \le C | b - a |$ for every $a, b \in [0,T]$ and $\mathcal{G}$--almost surely. Hence, \cite[Lemma 3.4]{trelat2000} and the dominated convergence theorem finally provide that
$$
\mathbb{E}\left[ \underset{t \in [0,T]}{\sup} \ \left\| \int^t_0 \left( u^i_k(s) - u^i(s) \right) b_i(s,x(s)) \; \mathrm{d}s \right\|^{\ell} \right] \rightarrow 0 , \quad k \rightarrow \infty .
$$
It is easy to conclude the proof by applying a routine Gr\"onwall inequality argument.
\end{proof}

The sought-after convergence of trajectories is a consequence of Lemma \ref{lemma_convTraj} when the conditions of Theorem \ref{theo_main} are met. In addition, limiting points for the sequences of controls fulfilling the conditions of Theorem \ref{theo_main} always exist up to some subsequence. Indeed, in this case, the set of admissible controls $\mathcal{U}$ is closed and convex for the strong topology of $L^2$. Hence it is closed for the weak topology of $L^2$. Therefore, since $(u_k)_{k \in \mathbb{N}}$ is bounded in $L^2$, there exists $u \in \mathcal{U}$ such that, up to some subsequence, $(u_k)_{k \in \mathbb{N}}$ weakly converges to $u$ for the weak topology of $L^2$. It remains to show that the process $x$ is feasible for OCP. For this, we compute
\begin{align*}
    \Big| \mathbb{E}[ &g(x(T))] \Big| \le \mathbb{E}\Big[ \| g(x(T)) - g(x_k(T)) \| \Big] + \mathbb{E}\left[ \left\| \frac{\partial g}{\partial x}(x_k(T)) (x_{k+1}(T) - x_k(T)) \right\| \right] \\
    &\le C \left( \mathbb{E}\Big[ \| x(T) - x_k(T) \|^2 \Big]^{\frac{1}{2}} + \mathbb{E}\Big[ \| x_{k+1}(T) - x(T) \|^2 \Big]^{\frac{1}{2}} \right) \rightarrow 0
\end{align*}
due to Lemma \ref{lemma_convTraj} and the dominated convergence theorem (here $C \ge 0$ is a constant coming from $(A_2)$, and we use the continuity of the sample paths of $x$).

\subsubsection{Convergence of Variational Inequalities} \label{sec_convVariational}

We start with a crucial result on linear stochastic differential equations. Recall the notation of Section \ref{sec_PMP}.
\begin{lmm}[Convergence of variational inequalities] \label{lemma_convZ}
Fix $\ell \ge 2$ and consider sequences of times $(r_k)_{k \in \mathbb{N}} \subseteq [0,T]$ and of uniformly bounded variables $(\tilde \xi_k)_{k \in \mathbb{N}}$ such that for every $k \in \mathbb{N}$, $\tilde \xi_k \in L^2_{\mathcal{F}_{r_k}}(\Omega;\mathbb{R}^{n+1})$. Assume that $r_k \rightarrow r$ with $r_k \le r$ for $k \in \mathbb{N}$ and $\tilde \xi_k \overset{L^{\ell}}{\rightarrow} \tilde \xi \in L^2_{\mathcal{F}_r}(\Omega;\mathbb{R}^{n+1})$ with $\tilde \xi$ bounded. Denote $\tilde w_{k+1}$, $\tilde w$ the stochastic process solutions, respectively, to
$$
\begin{cases}
    \displaystyle \mathrm{d}w(t) = \left( \frac{\partial b_0}{\partial x}(t,x_k(t)) + \sum^m_{i=1} u^i_k(t) \frac{\partial b_i}{\partial x}(t,x_k(t)) \right) w(t) \; \mathrm{d}t + \frac{\partial \sigma}{\partial x}(t,x_k(t)) w(t) \; \mathrm{d}B_t \\
    \displaystyle \mathrm{d}w^{n+1}(t) = \left( \frac{\partial H}{\partial x}(x_{k+1}(t)) + \frac{\partial L}{\partial x}(t,u_k(t),x_k(t)) \right) w(t) \; \mathrm{d}t \\
    \tilde w(t) = 0 , \ t \in [0,r_k) , \quad \tilde w(r_k) = \tilde \xi_{k+1} ,
\end{cases}
$$
$$
\begin{cases}
    \displaystyle \mathrm{d}w(t) = \left( \frac{\partial b_0}{\partial x}(t,x(t)) + \sum^m_{i=1} u^i(t) \frac{\partial b_i}{\partial x}(t,x(t)) \right) w(t) \; \mathrm{d}t + \frac{\partial \sigma}{\partial x}(t,x(t)) w(t) \; \mathrm{d}B_t \\
    \displaystyle \mathrm{d}w^{n+1}(t) = \left( \frac{\partial H}{\partial x}(x(t)) + \frac{\partial L}{\partial x}(t,u(t),x(t)) \right) w(t) \; \mathrm{d}t \\
    \tilde w(t) = 0 , \ t \in [0,r) , \quad \tilde w(r) = \tilde \xi .
\end{cases}
$$
Under the assumptions of Lemma \ref{lemma_convTraj}, it holds that $\mathbb{E}\left[ \underset{t \in [r,T]}{\sup} \ \| \tilde w_k(t) - \tilde w(t) \|^{\ell} \right] \rightarrow 0$ for $k \rightarrow \infty$.
\end{lmm}

\begin{proof}
From $r_k \le r$, $k \in \mathbb{N}$, for $t \in [r,T]$ we have (below, $C \ge 0$ is a constant)
\begin{align} \label{eq_proofConvLinear}
    \| &\tilde w_{k+1}(t) - \tilde w(t) \|^{\ell} \le C \| \tilde \xi_{k+1} - \tilde \xi \|^{\ell} \\
    & + C \left\| \int^t_{r_k} \left( \frac{\partial H}{\partial x}(x_{k+1}(s)) w_{k+1}(s) - \frac{\partial H}{\partial x}(x(s)) w(s) \right) \; \mathrm{d}s \right\|^{\ell} \nonumber \\
    & + C \left\| \int^t_{r_k} \left( \frac{\partial b_0}{\partial x}(s,x_k(s)) w_{k+1}(s) - \frac{\partial b_0}{\partial x}(s,x(s)) w(s) \right) \; \mathrm{d}s \right\|^{\ell} \nonumber \\
    & + C \left\| \int^t_0 \mathbbm{1}_{[r_k,T]}(t) \left( \frac{\partial \sigma}{\partial x}(s,x_k(s)) w_{k+1}(s) - \frac{\partial \sigma}{\partial x}(s,x(s)) w(s) \right) \; \mathrm{d}B_s \right\|^{\ell} \nonumber \\
    & + C \left\| \int^t_{r_k} \left( \frac{\partial L}{\partial x}(t,u_k(t),x_k(t)) w_{k+1}(s) - \frac{\partial L}{\partial x}(t,u(t),x(t)) w(s) \right) \; \mathrm{d}s \right\|^{\ell} \nonumber \\
    & + C \sum^{m}_{i=1} \left\| \int^t_{r_k} \left( u^i_k(s) \frac{\partial b_i}{\partial x}(s,x_k(s)) w_{k+1}(s) - u^i(s) \frac{\partial b_i}{\partial x}(s,x(s)) w(s) \right) \; \mathrm{d}s \right\|^{\ell} . \nonumber
\end{align}
We consider the expectation of the fourth term on the right-hand side (a similar argument can be developed for the first three terms, which are omitted in the interest of clarity and brevity). The Burkholder--Davis--Gundy and Hold\"er inequalities give
\begin{align*}
    &\mathbb{E}\bigg[ \bigg\| \int^t_0 \mathbbm{1}_{[r_k,T]}(t) \left( \frac{\partial \sigma}{\partial x}(s,x_k(s)) w_{k+1}(s) - \frac{\partial \sigma}{\partial x}(s,x(s)) w(s) \right) \; \mathrm{d}B_s \bigg\|^{\ell} \bigg] \le \\
    & \le C |r - r_k| + C \int^t_r \mathbb{E}\left[ \underset{s' \in [0,s]}{\sup} \ \left\| w_{k+1}(s') - w(s') \right\|^{\ell} \right] \; \mathrm{d}s + C \mathbb{E}\left[ \underset{s \in [0,T]}{\sup} \ \| w(s) \|^{2 \ell} \right] \mathbb{E}\left[ \underset{s \in [0,T]}{\sup} \ \| x_k(s) - x(s) \|^{2 \ell} \right] ,
\end{align*}
where the last inequality holds due to Lemma \ref{lemma_convTraj} and because Lemma \ref{lemma_bound} may be readily extended to $\tilde w$.

Finally, due to the fact that $L$ is affine with respect to the control variable, we may handle the fifth and the sixth terms in \eqref{eq_proofConvLinear} by combining the argument above with the final steps in the proof of Lemma \ref{lemma_convTraj}, and a routine Gr\"onwall inequality argument provides the conclusion.
\end{proof}

Consider $(u_k)_{k \in \mathbb{N}} \subseteq \mathcal{U}$ which converges to $u \in \mathcal{U}$ for the weak convergence of $L^2$. By assuming $(A_4)$, we may denote by $\mathcal{L} \subseteq [0,T]$ the full Lebesgue-measure subset such that $r \in \mathcal{L}$ if and only if $r$ is a Lebesgue point for $u$ and $r \notin \cup_{k \in \mathbb{N}} \mathcal{D}_{k+1}$ (we use the notation introduced with $(A_4)$). We prove the existence of a non-zero vector $\tilde{\mathfrak{p}} \in \mathbb{R}^{q+1}$ whose last component is non-positive such that for every $r \in \mathcal{L}$, $v \in U$,
\begin{equation} \label{eq_varIneqProof}
    \tilde{\mathfrak{p}}^{\top} \mathbb{E}\left[ \frac{\partial \tilde g}{\partial \tilde x}(\tilde x(T)) z_{r,v}(T) \right] \le 0 ,
\end{equation}
where the $\mathcal{F}$--adapted with continuous sample paths stochastic process $z_{r,v} : [0,T] \times \Omega \rightarrow \mathbb{R}^{n+1}$ solves \eqref{eq_rightSDE} with
$$
A(t) = \left( \begin{array}{c} \displaystyle \frac{\partial b_0}{\partial x}(t,x(t)) + \sum^m_{i = 1} u^i(t) \frac{\partial b_i}{\partial x}(t,x(t))  \\
\displaystyle \frac{\partial H}{\partial x}(x(t)) + \frac{\partial L}{\partial x}(t,u(t),x(t)) \end{array} \right) , \quad D(t) = \left( \begin{array}{c} \displaystyle \frac{\partial \sigma}{\partial x}(t,x(t))  \\
\displaystyle 0 \end{array} \right) ,
$$
$$
\tilde \xi_r = \tilde \xi_{r,v} \triangleq \tilde b(r,v,x(r)) - \tilde b(r,u(r),x(r)) ,
$$
where we denote $\tilde b = (b^{\top},f^0)^{\top}$ -- we will use this notation from now on.

\firstRev{For this, due to \eqref{eq_variationUpgraded}, for every $k \in \mathbb{N}$, the optimality of $(u_{k+1},x_{k+1})$ for LOCP$^{\Delta}_{k+1}$ provides a non-zero vector $(\tilde{\mathfrak{p}}^{\top}_{k+1},p^1_{k+1})^{\top} \in \mathbb{R}^{q+2}$ whose second to last component $p^0_{k+1} \le 0$ is non-zero by assumption, so that for $r \in \mathcal{L}$ and $v \in U$, it holds that
$$
\left(\begin{array}{c} \tilde{\mathfrak{p}}_{k+1} \\ p^1_{k+1} \end{array}\right)^{\top} \mathbb{E}\left[ \left(\begin{array}{cc} \displaystyle \frac{\partial \tilde g}{\partial \tilde x}(\tilde x_k(T)) & 0 \\ 0 & 1 \end{array}\right) \left(\begin{array}{c} z^{k+1}_{r,v}(T) \\ (z^{k+1}_{r,v})^{n+2}(T) \end{array}\right) \right] \le 0 ,
$$
where with the notation
$$
A_{k+1}(t) = \left( \begin{array}{cc} \displaystyle \frac{\partial b_0}{\partial x}(t,x_k(t)) + \sum^m_{i = 1} u^i_k(t) \frac{\partial b_i}{\partial x}(t,x_k(t)) & 0 \\
\displaystyle \frac{\partial H}{\partial x}(x_{k+1}(t)) + \frac{\partial L}{\partial x}(t,u_k(t),x_k(t)) & 0 \\
2 ( x_{k+1}(t) - x_k(t) ) & 0 \end{array} \right) , \quad D_{k+1}(t) = \left( \begin{array}{cc} \displaystyle \frac{\partial \sigma}{\partial x}(t,x_k(t)) & 0 \\
\displaystyle 0 & 0 \\
\displaystyle 0 & 0 \end{array} \right) ,
$$
the $\mathcal{F}$--adapted with continuous sample paths stochastic process $((z^k_{r,v})^{\top},(z^k_{r,v})^{n+2})^{\top} : [0,T] \times \Omega \rightarrow \mathbb{R}^{n+2}$ solves a higher dimensional version of \eqref{eq_rightSDE} with
$$
A = A_{k+1} , \ D = D_{k+1} , \ \textnormal{and initial condition} \ \left(\begin{array}{c} \tilde \xi^{k+1}_{r,v} \\ 0 \end{array}\right) \triangleq \left(\begin{array}{c} \tilde b_{k+1}(r,v,x_{k+1}(r)) - \tilde b_{k+1}(r,u_{k+1}(r),x_{k+1}(r)) \\ 0 \end{array}\right) ,
$$
where again we denote $\tilde b_{k+1} = (b^{\top}_{k+1},f^0_{k+1})^{\top}$.}

Now, fix $r \in \mathcal{L}$ and $v \in U$. The following comes from \firstRev{combining Lemma \ref{lemma_continuity} with \cite[Lemma 3.11]{bonalli2019}}.
\begin{lmm}[Pointwise convergence on controls] \label{lemma_convContr}
Under $(A_4)$, there exists $(r_k)_{k \in \mathbb{N}} \subseteq (0,r)$ such that for every $k \in \mathbb{N}$, $r_k$ is a Lebesgue point for $u_k$, and $r_k \rightarrow r$, $u_k(r_k) \rightarrow u(r)$ as $k \rightarrow \infty$.
\end{lmm}
If $(r_k)_{k \in \mathbb{N}} \subseteq (0,r)$ denotes the sequence given by Lemma \ref{lemma_convContr}, we define $\tilde \xi_k = \tilde \xi^k_{r_k,v}$ and $\tilde \xi = \tilde \xi_{r,v}$. Straightforward computations give (below, $C \ge 0$ is a constant)
\begin{align} \label{eq_boundProof}
    \mathbb{E}\left[ \| \tilde \xi_{k+1} - \tilde \xi \|^2 \right] \le &C \bigg( \| u_{k+1}(r_{k+1}) - u(r) \|^2 + | r_{k+1} - r |^2 \\
    & + \mathbb{E}\left[ \| x(r_{k+1}) - x(r) \|^2 \right] + \mathbb{E}\left[ \underset{s \in [0,T]}{\sup} \ \| x_k(s) - x(s) \|^2 \right] \bigg) , \nonumber
\end{align}
from which $\mathbb{E}\left[ \| \tilde \xi_{k+1} - \tilde \xi \|^2 \right] \rightarrow 0$ for $k \rightarrow \infty$ (Lemma \ref{lemma_convTraj} and \ref{lemma_convContr}). Therefore, from Lemma \ref{lemma_convZ} we infer that
\firstRev{$$
\mathbb{E}\left[ \underset{t \in [r,T]}{\sup} \ \| z^k_{r_k,v}(t) - z_{r,v}(t) \|^2 \right] \rightarrow 0 , \quad k \rightarrow \infty ,
$$
which, together with $\Delta_k \rightarrow 0$, readily yields
$$
\mathbb{E}\left[ | (z^k_{r_k,v})^{n+2}(T) | \right] \rightarrow 0 , \quad k \rightarrow \infty .
$$
At this step, we point out that the variational inequalities in \eqref{eq_variation} and in \eqref{eq_variationUpgraded} still hold if we take multipliers of norm one. Specifically, we may assume that $\| (\tilde{\mathfrak{p}}^{\top}_{k+1},p^1_{k+1})^{\top} \| = 1$ for every $k \in \mathbb{N}$. Therefore, up to some subsequence, there exists a vector $(\tilde{\mathfrak{p}}^{\top},p^1)^{\top} = (\mathfrak{p}^{\top},p^0,p^1)^{\top} \in \mathbb{R}^{q+1}$ such that $(\tilde{\mathfrak{p}}^{\top}_k,p^1_k)^{\top} \rightarrow (\tilde{\mathfrak{p}}^{\top},p^1)^{\top}$ for $k \rightarrow \infty$ and satisfying $(\tilde{\mathfrak{p}}^{\top},p^1)^{\top} \neq 0$, $p^0 \le 0$. We use this remark to conclude as follows. The definition of $\tilde g$ and the H\"older inequality give ($C \ge 0$ is a constant)
\begin{align*}
    &\tilde{\mathfrak{p}}^{\top} \mathbb{E}\left[ \frac{\partial \tilde g}{\partial \tilde x}(\tilde x(T)) z_{r,v}(T) \right] \le \left| \tilde{\mathfrak{p}}^{\top} \mathbb{E}\left[ \frac{\partial \tilde g}{\partial \tilde x}(\tilde x(T)) z_{r,v}(T) \right] - \left(\begin{array}{c} \tilde{\mathfrak{p}}_{k+1} \\ p^1_{k+1} \end{array}\right)^{\top} \mathbb{E}\left[ \left(\begin{array}{cc} \displaystyle \frac{\partial \tilde g}{\partial \tilde x}(\tilde x_k(T)) & 0 \\ 0 & 1 \end{array}\right) \left(\begin{array}{c} z^{k+1}_{r_{k+1},v}(T) \\ (z^{k+1}_{r_{k+1},v})^{n+2}(T) \end{array}\right) \right] \right| \\
    & \le C \bigg( \| \tilde{\mathfrak{p}} - \tilde{\mathfrak{p}}_{k+1} \| + p^1_{k+1} \mathbb{E}\left[ | (z^k_{r_k,v})^{n+2}(T) | \right]  \\
    &\hspace{30ex} + \mathbb{E}\left[ \| z^{k+1}_{r_{k+1},v}(T) - z_{r,v}(T) \|^2 \right]^{\frac{1}{2}} + \mathbb{E}\Big[ \| z_{r,v}(T) \|^2 \Big]^{\frac{1}{2}} \mathbb{E}\left[ \underset{s \in [0,T]}{\sup} \ \| x_k(s) - x(s) \|^2 \right]^{\frac{1}{2}} \bigg) ,
\end{align*}
and in this case, \eqref{eq_varIneqProof} follows from Lemma \ref{lemma_convTraj} and the convergences obtained above. From what we showed in Section \ref{sec_daulDerivation}, this latter inequality yields that $(u,\mathfrak{p},p^0)$ is extremal for OCP as soon as $(\mathfrak{p}^{\top},p^0)^{\top} \neq 0$.}

Finally, we turn to the case for which the sequence $(u_k)_{k \in \mathbb{N}} \subseteq \mathcal{U}$ converges to $u \in \mathcal{U}$ for the strong topology of $L^2$, but without assuming $(A_4)$. For this, fix $v \in U$ and define the stochastic processes $\tilde \xi_k(s) = \tilde \xi^k_{s,v}$ and $\tilde \xi(s) = \tilde \xi_{s,v}$, where $s \in [0,T]$. Similar computations to the ones developed to compute the bound \eqref{eq_boundProof} provide that $\int^T_0 \mathbb{E}\left[ \| \tilde \xi_{k}(s) - \tilde \xi(s) \|^2 \right] \; \mathrm{d}s \rightarrow 0$ as $k \rightarrow \infty$, and therefore up to some subsequence, the quantity $\mathbb{E}\left[ \| \tilde \xi_{k}(s) - \tilde \xi(s) \|^2 \right]$ converges to zero for $k \rightarrow 0$ a.e. in $[0,T]$. By taking countable intersections of sets of Lebesgue points (one for each control $u_k$, for all $k \in \mathbb{N}$), it follows that the argument above can be iterated exactly in the same manner (via Lemma \ref{lemma_convZ}), leading to the same conclusion.

\subsubsection{Convergence of Multipliers and Conclusion}

By applying the construction developed in Section \ref{sec_daulDerivation} to the variational inequality \eqref{eq_varIneqProof}, we retrieve a tuple \firstRev{$(p,\mathfrak{p},p^0,q)$, where $\mathfrak{p} \in \mathbb{R}^q$, $p^0 \le 0$ are constant, $p \in L^2_{\mathcal{F}}(\Omega;C^0([0,T];\mathbb{R}^n))$ and $q \in L^2_{\mathcal{F}}([0,T] \times \Omega;\mathbb{R}^n)$, such that, as soon as $(\mathfrak{p}^{\top},p^0)^{\top} \neq 0$, the tuple $(u,\mathfrak{p},p^0,q)$ is extremal for OCP associated with $(u,x,p,\mathfrak{p},p^0,q)$ satisfying conditions 1., 2., and 3. of Theorem \ref{Theo_PMP}}. To conclude, it only remains to prove that
\begin{equation} \label{eq_lastConv}
\mathbb{E}\left[ \underset{s \in [0,T]}{\sup} \ \| p_k(s) - p(s) \|^2 \right] \rightarrow 0 , \quad \int^T_0 \mathbb{E}\Big[ \| q_k(s) - q(s) \|^2 \Big] \; \mathrm{d}s \rightarrow 0 , \quad k \rightarrow \infty.
\end{equation}
Here, each tuple \firstRev{$(u_k,\mathfrak{p}_k,p^0_k)$ is the extremal of LOCP$^{\Delta}_k$ associated with the tuple $(u_k,x_k,p^k,\mathfrak{p}_k,p^0_k,q_k)$, where $(p_k,q_k)$ solves the adjoint equation of Theorem \ref{Theo_PMPk}. Since, under the assumption $(\mathfrak{p}^{\top},p^0)^{\top} \neq 0$, the multipliers $p^1_k$, $k \in \mathbb{N}$, and $p^1$ do not play any role, for the sake of clarity and brevity of notation, and without loss of generality, in what follows we implicitly assume $p^1_k = p^1 = 0$, for every $k \in \mathbb{N}$.}

Let us start with the first convergence of \eqref{eq_lastConv}. For this, fix $k \in \mathbb{N}$ and consider the process
$$
\left( \mathbb{E}\left[ \underset{s \in [0,T]}{\sup} \ \left\| \frac{\partial \tilde g}{\partial \tilde x}(\tilde x_k(T)) \phi_{k+1}(T) \psi_{k+1}(s) - \frac{\partial \tilde g}{\partial \tilde x}(\tilde x(T)) \phi(T) \psi(s) \right\| \bigg| \mathcal{F}_t \right] \right)_{t \in [0,T]} ,
$$
where $\phi_{k+1}$, $\psi_{k+1}$ solve \eqref{eq_phiPsi} with matrices $A_{k+1}$, $D_{k+1}$ \firstRev{from which, due to $p^1_k = p^1 = 0$, for every $k \in \mathbb{N}$, we remove the last column and row}, whereas $\phi$, $\psi$ solve \eqref{eq_phiPsi} with matrices $A$, $D$, those matrices being defined above. Due to a straightforward extension of Lemma \ref{lemma_bound} to equations \eqref{eq_phiPsi}, this process is a martingale, bounded in $L^2$. Hence, the martingale representation theorem allows us to infer that this process is a martingale with continuous sample paths, and Doob and Jensen inequalities give
{
\begin{align} \label{eq_martingaleProof}
    \mathbb{E}&\left[ \underset{t \in [0,T]}{\sup} \ \mathbb{E}\left[ \underset{s \in [0,T]}{\sup} \ \left\| \frac{\partial \tilde g}{\partial \tilde x}(\tilde x_k(T)) \phi_{k+1}(T) \psi_{k+1}(s) - \frac{\partial \tilde g}{\partial \tilde x}(\tilde x(T)) \phi(T) \psi(s) \right\| \bigg| \mathcal{F}_t \right]^2 \right] \le \\
    & \le 4 \mathbb{E}\left[ \underset{s \in [0,T]}{\sup} \ \left\| \frac{\partial \tilde g}{\partial \tilde x}(\tilde x_k(T)) \phi_{k+1}(T) \psi_{k+1}(s) - \frac{\partial \tilde g}{\partial \tilde x}(\tilde x(T)) \phi(T) \psi(s) \right\|^2 \right] , \nonumber
\end{align}}
which holds for every $k \in \mathbb{N}$. By combining \eqref{eq_martingaleProof} with \eqref{eq_adjointVector}, we compute
\begin{align*}
    \mathbb{E}&\left[ \underset{s \in [0,T]}{\sup} \ \| p_{k+1}(s) - p(s) \|^2 \right] \le C \mathbb{E}\left[ \underset{t \in [0,T]}{\sup} \ \left\| \frac{\partial \tilde g}{\partial \tilde x}(\tilde x(t_f)) \phi(T) \psi(s) \right\|^2 \right] \| \tilde{\mathfrak{p}}_{k+1} - \tilde{\mathfrak{p}} \|^2 \\
    & + C \mathbb{E}\left[ \underset{s \in [0,T]}{\sup} \ \left\| \frac{\partial \tilde g}{\partial \tilde x}(\tilde x_k(T)) \phi_{k+1}(T) \psi_{k+1}(s) - \frac{\partial \tilde g}{\partial \tilde x}(\tilde x(T)) \phi(T) \psi(s) \right\|^2 \right] ,
\end{align*}
where $C \ge 0$ is a constant. Up to some subsequence, the first term on the right-hand side converges to zero. Moreover, the definition of $\tilde g$ and H\"older inequality give
{\small
\begin{align*}
    &\mathbb{E}\left[ \underset{s \in [0,T]}{\sup} \ \left\| \frac{\partial \tilde g}{\partial \tilde x}(\tilde x_k(T)) \phi_{k+1}(T) \psi_{k+1}(s) - \frac{\partial \tilde g}{\partial \tilde x}(\tilde x(T)) \phi(T) \psi(s) \right\|^2 \right] \le \\
    &\le C \mathbb{E}\left[ \| \phi_{k+1}(T) \|^4 \right]^{\frac{1}{2}} \mathbb{E}\left[ \underset{s \in [0,T]}{\sup} \ \| \psi_{k+1}(s) - \psi(s) \|^4 \right]^{\frac{1}{2}} \\
    & \ + C \left( \mathbb{E}\left[ \underset{s \in [0,T]}{\sup} \ \| \psi(s) \|^4 \right]^{\frac{1}{2}} \mathbb{E}\left[ \underset{s \in [0,T]}{\sup} \ \| \phi_{k+1}(s) - \phi(s) \|^4 \right]^{\frac{1}{2}} + \mathbb{E}\left[ \underset{s \in [0,T]}{\sup} \ \| \phi(T) \psi(s) \|^4 \right]^{\frac{1}{2}} \mathbb{E}\left[ \underset{s \in [0,T]}{\sup} \ \| x_k(s) - x(s) \|^4 \right]^{\frac{1}{2}} \right) ,
\end{align*}}
and a straightforward extension of Lemma \ref{lemma_convZ} to equations \eqref{eq_phiPsi} entails that all the terms on the right-hand side tend to zero. The convergence of $(p_k)_{k \in \mathbb{N}}$ is proved.

It remains to prove the second convergence of \eqref{eq_lastConv}. For this, we apply It\^o formula to $\| p_{k+1}(t) - p(t) \|^2$, which due to \eqref{eq_adjointVector}, \eqref{eq_martingaleProof} and Lemma \ref{lemma_bound} extended to \eqref{eq_phiPsi} gives (below, $C \ge 0$ is a constant)
{
\begin{align*}
    &\mathbb{E}\Big[ \| p_{k+1}(0) - p(0) \|^2 \Big] + \int^T_0 \mathbb{E}\Big[ \| q_{k+1}(s) - q(s) \|^2 \Big] \; \mathrm{d}s = \mathbb{E}\Big[ \| p_{k+1}(T) - p(T) \|^2 \Big] \\
    &+ 2 \mathbb{E}\bigg[ \int^T_0 \Big( p_{k+1}(s) - p(s) \Big)^{\top} \bigg( p(s)^{\top}_{k+1} \frac{\partial b_{k+1}}{\partial x}(s,u_{k+1}(s),x_{k+1}(s)) - p(s)^{\top} \frac{\partial b}{\partial x}(s,u(s),x(s)) \bigg) \; \mathrm{d}s \bigg] \\
    &+ 2 \mathbb{E}\left[ \int^T_0 \Big( p_{k+1}(s) - p(s) \Big)^{\top} \bigg( p^0_{k+1} \frac{\partial f^0_{k+1}}{\partial x}(s,u_{k+1}(s),x_{k+1}(s)) - p^0 \frac{\partial f^0}{\partial x}(s,u(s),x(s)) \bigg) \; \mathrm{d}s \right] \\
    &+ 2 \mathbb{E}\left[ \int^T_0 \Big( p_{k+1}(s) - p(s) \Big)^{\top} \bigg( q(s)^{\top}_{k+1} \frac{\partial \sigma_{k+1}}{\partial x}(s,x_{k+1}(s)) - q(s)^{\top} \frac{\partial \sigma}{\partial x}(s,x(s)) \bigg) \; \mathrm{d}s \right] \\
    &\le C \mathbb{E}\left[ \underset{s \in [0,T]}{\sup} \ \| p_{k+1}(s) - p(s) \|^2 \right]^{\frac{1}{2}} \left( \int^T_0 \mathbb{E}\Big[ \| q_{k+1}(s) - q(s) \|^2 \Big] \; \mathrm{d}s \right)^{\frac{1}{2}} \\
    & \quad + C \mathbb{E}\left[ \underset{s \in [0,T]}{\sup} \ \| p_{k+1}(s) - p(s) \|^2 \right]^{\frac{1}{2}}\left( 1 + \left( \int^T_0 \mathbb{E}\Big[ \| q(s) \|^2 \Big] \; \mathrm{d}s \right)^{\frac{1}{2}} \right) .
\end{align*}}
The conclusion finally follows from Young's inequality and the convergence of $(p_k)_{k \in \mathbb{N}}$.

\section{Extension to Problems with Free Final Time and Stochastic Controls} \label{sec_extension}

The \firstRev{accumulation} properties of SCP can be extended to OCPs with free final time and whereby optimization is undertaken over stochastic controls. Specifically, in this section we investigate how results from Theorem \ref{theo_main} may be extended to finite-horizon, finite-dimensional non-linear stochastic General Optimal Control Problems (GOCP) having control-affine dynamics and uncontrolled diffusion, of the form
{%
\begin{equation*}
    \begin{cases}
        \displaystyle \underset{u , t_f}{\min} \ \mathbb{E}\left[ \int^{t_f}_0 f^0(s,u(s),x(s)) \; \mathrm{d}s \right] \triangleq \mathbb{E}\left[ \int^{t_f}_0 \Big( G(u(s)) + H(x(s)) + L(s,u(s),x(s)) \Big) \; \mathrm{d}s \right] \medskip \\
        \mathrm{d}x(t) = b(t,u(t),x(t)) \; \mathrm{d}t + \sigma(t,x(t)) \; \mathrm{d}B_t , \quad x(0) = x^0 , \quad \mathbb{E}\left[ g(x(t_f)) \right] = 0 ,
    \end{cases}
\end{equation*}}
where the final time $0 \le t_f \le T$ may be free or not (here, $T > 0$ is some fixed maximal time), and we optimize over controls $u \in \mathcal{U}$ which are either deterministic, i.e., $\mathcal{U} = L^2([0,T];U)$, or stochastic, i.e., $\mathcal{U} = L^2_{\mathcal{F}}([0,T] \times \Omega;U)$. Solutions to GOCP will be denoted by $(t_f,u,x) \in [0,T] \times \mathcal{U} \times L^{\ell}_{\mathcal{F}}(\Omega;C([0,T];\mathbb{R}^n))$, for $\ell \ge 2$.

By adopting the notation given in \eqref{eq_LCost}, the stochastic Linearized General Optimal Control Problem (LGOCP$^{\Delta}_{k+1}$) at iteration $k \in \mathbb{N}$ may be defined accordingly as
\begin{equation*}
    \begin{cases}
        \displaystyle \underset{u , t_f}{\min} \ \mathbb{E}\left[ \int^{t_f}_0 f^0_{k+1}(s,u(s),x(s)) \; \mathrm{d}s \right] \triangleq \mathbb{E}\bigg[ \int^{t_f}_0 f^0_{u_k,x_k}(s,u(s),x(s)) \; \mathrm{d}s \bigg] \medskip \\
        \mathrm{d}x(t) = b_{k+1}(t,u(t),x(t)) \; \mathrm{d}t + \sigma_{k+1}(t,x(t)) \; \mathrm{d}B_t , \quad x(0) = x^0 \medskip \\
        \hspace{6ex} \triangleq b_{u_k,x_k}(t,u(t),x(t)) \; \mathrm{d}t + \sigma_{x_k}(t,x(t)) \; \mathrm{d}B_t \medskip \\
        \displaystyle \mathbb{E}\left[ g_{k+1}(x(t_f)) \right] \triangleq \mathbb{E}\left[ g(x_k(t^k_f)) + \frac{\partial g}{\partial x}(x_k(t^k_f)) (x(t_f) - x_k(t^k_f)) \right] = 0 \medskip \\
        \displaystyle \int^T_0 \mathbb{E}\left[ \| x(s) - x_k(s) \|^2 \right] \; \mathrm{d}s \le \Delta_{k+1} , \quad |t_f - t^k_f| \le \Delta_{k+1} ,
    \end{cases}
\end{equation*}
whose solutions are denoted by $(t^{k+1}_f,u_{k+1},x_{k+1}) \in [0,T] \times \mathcal{U} \times L^{\ell}_{\mathcal{F}}(\Omega;C([0,T];\mathbb{R}^n))$, for $\ell \ge 2$. In particular, trust-region constraints are now imposed on the variable $t_f$ as well to derive convergence guarantees. 

We point out that the results we present in this section are not to be considered as part of the main contribution, but they rather aim at providing insights for the future design of efficient numerical methods for stochastic optimal control problems which consider free final time and stochastic admissible controls. In particular, extending SCP in the presence of free final time and stochastic controls requires introducing additional assumptions which might seem demanding. We leave the investigation of the validity of this framework under sharper assumptions as a future direction of research.

\subsection{Refined assumptions and extended result of convergence}

The presence of the free final time hinders the well-posedness of the stochastic PMP when applied to each LGOCP$^{\Delta}_{k+1}$, i.e., Theorem \ref{Theo_PMP}, and in turn the validity of Theorem \ref{theo_main} in such more general setting. To overcome this issue, we need to tighten assumptions \firstRev{$(A_1)$--$(A_3)$}. Specifically, although $(A_1)$ remains unchanged, assumptions \firstRev{$(A_2)$--$(A_3)$} are replaced by the followings, respectively:

\vspace{5pt}

\noindent $(A'_2)$ Mappings $g$, $H$, and $L_i$, $i=0,\dots,m$, either are affine-in-state or have compact supports in $\mathbb{R}^n$ and in $[0,T] \times \mathbb{R}^n$, respectively. In the case of free final time, $g$ is affine.

\vspace{5pt}

\noindent $(A'_3)$ For every $k \in \mathbb{N}$, LGOCP$^{\Delta}_{k+1}$ is feasible. In addition, in the case of free final time, for every $k \in \mathbb{N}$, any optimal control $u_k$ for LGOCP$^{\Delta}_k$ has continuous sample paths at the optimal final time $t^{k+1}_f$ for LGOCP$^{\Delta}_{k+1}$.

\vspace{5pt}

\firstRev{Assumptions $(A'_3)$ plays a crucial role to provide the existence of pointwise necessary conditions for optimality for the linearized problems LGOCP$^{\Delta}_k$, thus enabling classical formulations of the stochastic PMP. However, we do recognize $(A'_3)$ might be demanding and, as future research direction, we propose to relax it by leveraging integral-type necessary conditions for optimality, which are best suited to deal with optimal control problems which show explicit discontinuous dependence on time (e.g., \cite{bourdin2013}). Unfortunately, when optimization over stochastic controls is adopted, successfully leveraging weakly converging subsequences of controls to show ``weak'' guarantee of success for SCP, \firstRev{i.e., the equivalent of Theorem \ref{theo_main2}}, becomes more challenging, in which case only Theorem \ref{theo_main} still holds\firstRev{, though under appropriate modifications}. We leave proving success guarantees for SCP when optimizing over stochastic controls with weaker assumptions as a future direction of research.}

\firstRev{Extending the stochastic PMP to GOCP and LGOCP$^{\Delta}_{k+1}$ comes by assuming $(A_1)$, $(A'_2)$, and $(A'_3)$. In particular, under those assumptions we have the following extension of Theorem \ref{Theo_PMP} and of Theorem \ref{Theo_PMPk}}:
\begin{thrm}[Stochastic Pontryagin Maximum Principle for GOCP] \label{Theo_GPMP}
Let $(t_f,u,x)$ be a locally-optimal solution to GOCP. There exist $p\in L^2_{\mathcal{F}}(\Omega;C([0,T];\mathbb{R}^n))$ and a tuple $(\mathfrak{p},p^0,q)$, where $\mathfrak{p} \in \mathbb{R}^q$ and $p^0 \le 0$ are constant, and $q = (q_1,\dots,q_d) \in L^2_{\mathcal{F}}([0,T] \times \Omega;\mathbb{R}^{n \times d})$ such that the following relations are satisfied:
\begin{enumerate}
    \item Non-Triviality Condition: $(\mathfrak{p},p^0) \neq 0$.
    \item Adjoint Equation:
    \firstRev{\begin{align*}
        \mathrm{d}p(t) &= \displaystyle -\frac{\partial H}{\partial x}(t,u(t),x(t),p(t),p^0,q(t)) \; \mathrm{d}t + q(t) \; \mathrm{d}B_t , \quad p(t_f) = \displaystyle \mathbb{E}\left[ \frac{\partial g}{\partial x}(x(t_f)) \right]^{\top} \mathfrak{p} \in \mathbb{R}^n .
    \end{align*}}
    \item Maximality Condition:
    $$
    \firstRev{u(t) = \underset{v \in U}{\arg \max} \ \mathbb{E}\Big[ H(t,v,x(t),p(t),p^0,q(t)) \Big]} , \ \textnormal{a.e.} \ \textnormal{(deterministic controls)}
    $$
    $$
    \firstRev{u(t) = \underset{v \in U}{\arg \max} \ H(t,v,x(t),p(t),p^0,q(t))} , \ \textnormal{a.e., \ a.s.} \ \textnormal{(stochastic controls)}
    $$
    \item Transversality Condition: if the final time is free
    $$
    \firstRev{\underset{v \in U}{\max} \ \mathbb{E}\Big[ H(t_f,v,x(t_f),p(t_f),p^0,q(t_f)) \Big] \ge 0} \ \textnormal{(deterministic controls)}
    $$
    $$
    \firstRev{\mathbb{E}\left[ \underset{v \in U}{\max} \ H(t_f,v,x(t_f),p(t_f),p^0,q(t_f)) \right] \ge 0} \ \textnormal{(stochastic controls)}
    $$
    where equalities hold in the case $t_f < T$.
\end{enumerate}
\firstRev{The quantity $(t_f,u,\mathfrak{p},p^0)$ uniquely determines $x$, $p$, and $q$ and is called extremal for GOCP (associated with the tuple $(t_f,u,x,p,\mathfrak{p},p^0,q)$, or simply with $(t_f,u,x)$). An extremal $(t_f,u,\mathfrak{p},p^0)$ is called \textnormal{normal} if $p^0 \neq 0$}.
\end{thrm}

\begin{thrm}[Weak Stochastic Pontryagin Maximum Principle for LGOCP$^{\Delta}_k$] \label{Theo_GPMPk}
\firstRev{Let $(t^k_f,u_k,x_k)$ be a locally-optimal solution to LGOCP$^{\Delta}_k$. There exist $p_k \in L^2_{\mathcal{F}}(\Omega;C([0,T];\mathbb{R}^n))$, a tuple $(\mathfrak{p}_k,p^0_k,p^1_k)$, where $\mathfrak{p}_k \in \mathbb{R}^q$, $p^0_k \le 0$, and $p^1_k \in \mathbb{R}$ are constant, and $q_k = ((q_k)_1,\dots,(q_k)_d) \in L^2_{\mathcal{F}}([0,T] \times \Omega;\mathbb{R}^{n \times d})$ such that the following relations are satisfied:
\begin{enumerate}
    \item Non-Triviality Condition: $(\mathfrak{p}_k,p^0_k,p^1_k) \neq 0$.
    \item Adjoint Equation:
    \begin{align*}
        \mathrm{d}p_k(t) &= \displaystyle -\frac{\partial H_k}{\partial x}(t,u_k(t),x_k(t),p_k(t),p^0_k,p^1_k,q_k(t)) \; \mathrm{d}t + q_k(t) \; \mathrm{d}B_t , \quad p_k(t^k_f) = \displaystyle \mathbb{E}\left[ \frac{\partial g_k}{\partial x}(x_k(t^k_f)) \right]^{\top} \mathfrak{p}_k \in \mathbb{R}^n .
    \end{align*}
    \item Maximality Condition:
    $$
    u_k(t) = \underset{v \in U}{\arg \max} \ \mathbb{E}\Big[ H_k(t,v,x_k(t),p_k(t),p^0_k,p^1_k,q_k(t)) \Big] , \ \textnormal{a.e.} \ \textnormal{(deterministic controls)}
    $$
    $$
    u_k(t) = \underset{v \in U}{\arg \max} \ H_k(t,v,x_k(t),p_k(t),p^0_k,p^1_k,q_k(t)) , \ \textnormal{a.e., \ a.s.} \ \textnormal{(stochastic controls)}
    $$
    \item Transversality Condition: if the final time is free
    $$
    \underset{v \in U}{\max} \ \mathbb{E}\Big[ H_k(t^k_f,v,x_k(t^k_f),p_k(t^k_f),p^0_k,p^1_k,q_k(t^k_f)) \Big] \ge 0 \ \textnormal{(deterministic controls)}
    $$
    $$
    \mathbb{E}\left[ \underset{v \in U}{\max} \ H_k(t^k_f,v,x_k(t^k_f),p_k(t^k_f),p^0_k,p^1_k,q_k(t^k_f)) \right] \ge 0 \ \textnormal{(stochastic controls)}
    $$
    where equalities hold in the case $t^k_f < T$.
\end{enumerate}
The quantity $(t^k_f,u_k,\mathfrak{p}_k,p^0_k,p^1_k)$ uniquely determines $x_k$, $p_k$, and $q_k$ and is called extremal for LOCP$^{\Delta}_k$ (associated with $(t^k_f,u_k,x_k,p_k,\mathfrak{p}_k,p^0_k,p^1_k,q_k)$, or $(t^k_f,u_k,x_k)$). An extremal $(t^k_f,u_k,\mathfrak{p}_k,p^0_k,p^1_k)$ is called \textnormal{normal} if $p^0_k \neq 0$}.
\end{thrm}

Extending the stochastic PMP to the new linearized problems LGOCP$^{\Delta}_k$, and in particular the new transversality condition 4., additionally requires to assume $(A'_3)$ together with $(A_1)$ and $(A'_2)$. Although the proof of this result for fixed-final-time problems is well-established (see \cite[Chapter 3]{Yong1999}), we could not find any published proof of Theorem \ref{Theo_GPMP} when the final time is free. Therefore, we provide its proof in Section \ref{sec_proofGPMP}. \firstRev{The proof of Theorem \ref{Theo_GPMPk} is achieved similarly, thus for the sake of clarity and brevity, we avoid reporting any detail concerning the latter}. Thanks to Theorem \ref{Theo_GPMP}, we can extend the convergence of SCP as follows:

\begin{thrm}[Generalized Properties of Accumulation Points for SCP] \label{theo_mainGeneralized}
\firstRev{Assume that $(A_1)$, $(A'_2)$, and $(A'_3)$ hold and that SCP generates a sequence $(\Delta_k,t^k_f,u_k,x_k)_{k \in \mathbb{N}}$ such that $(\Delta_k)_{k \in \mathbb{N}} \subseteq \mathbb{R}_+ \setminus \{ 0 \}$ converges to zero, and for every $k \ge 1$, the tuple $(t^k_f,u_k,x_k)$ locally solves LGOCP$^{\Delta}_k$. For every $k \ge 1$, letting $(t^k_f,u_k,\mathfrak{p}_k,p^0_k,p^1_k)$ be an extremal associated with $(t^k_f,u_k,x_k)$ for LGOCP$^{\Delta}_k$ (whose existence is ensured by Theorem \ref{Theo_GPMPk}), assume the following Accumulation Condition holds:
\begin{itemize}
    \item[\textnormal{(AC)}] Up to some subsequence, $(t^k_f,u_k,\mathfrak{p}_k,p^0_k,p^1_k)$ converges to some $(u,\mathfrak{p},p^0,p^1) \in L^2([0,T];\mathbb{R}^m) \times \mathbb{R}^{q+2}$ for the strong topology of $L^2_{\mathcal{F}}([0,T]\times\Omega;\mathbb{R}^m) \times \mathbb{R}^{q+2}$.
\end{itemize}
If $(\mathfrak{p},p^0) \neq 0$, then $(t^k_f,u,\mathfrak{p},p^0)$ is an extremal for OCP associated with $(t^k_f,u,x_u)$.}
\end{thrm}

The guarantees offered by Theorem \ref{theo_mainGeneralized} read similarly to what we have explained in Section \ref{sec_mainDiscussion}. One sees that the computations provided in Section \ref{sec_proofConv} for the proof of Theorem \ref{theo_main} straightforwardly generalize \firstRev{as soon as weak convergence of controls in $L^2([0,T];\mathbb{R}^m)$ is replaced with strong convergence of controls in $L^2_{\mathcal{F}}([0,T]\times\Omega;\mathbb{R}^m)$ (actually, the proofs become even simpler),} and may be endorsed to prove Theorem \ref{theo_mainGeneralized}, provided that we are capable to extend the proof of Theorem \ref{Theo_PMP} \firstRev{and Theorem \ref{Theo_PMPk}} in the case of free final time and stochastic controls. Therefore, to conclude, in the next section we develop the necessary technical details which enable proving Theorem \ref{Theo_GPMP} \firstRev{and Theorem \ref{Theo_GPMPk}} through the machinery developed in Section \ref{sec_proofConv}. \firstRev{Specifically, since the proofs of Theorem \ref{Theo_GPMP} and Theorem \ref{Theo_GPMPk} are similar, for the sake of clarity and brevity we only provide details for the proof of Theorem \ref{Theo_GPMP}.}

\subsection{Proof of the extension for the stochastic Pontryagin Maximum Principle} \label{sec_proofGPMP}

Before getting started, we need to introduce the notion of a \textit{Lebesgue point} for a stochastic control $u \in L^2_{\mathcal{F}}([0,T] \times \Omega;\mathbb{R}^m)$. For this, we adopt the theory of Bochner integrals, showing that $L^2_{\mathcal{F}}([0,T] \times \Omega;\mathbb{R}^m) \subseteq L^2([0,T];L^2_{\mathcal{G}}(\Omega;\mathbb{R}^m))$, where the latter is the space of Bochner integrable mappings $u : [0,T] \rightarrow L^2_{\mathcal{G}}(\Omega;\mathbb{R}^m)$. Let $u \in L^2_{\mathcal{F}}([0,T] \times \Omega;\mathbb{R}^m)$. First of all, by definition we see that for almost every $t \in [0,T]$, it holds that $\mathbb{E}\left[ \| u(t) \|^2 \right] < \infty$, and thus this control is well-defined as a mapping $u : [0,T] \rightarrow L^2_{\mathcal{G}}(\Omega;\mathbb{R}^m)$. Since $\Omega$ is second-countable, $L^2_{\mathcal{G}}(\Omega;\mathbb{R}^m)$ is separable, and therefore the claim follows from the Pettis measurability theorem once we prove that $u : [0,T] \rightarrow L^2_{\mathcal{G}}(\Omega;\mathbb{R}^m)$ is strongly measurable with respect to the Lebesgue measure of $[0,T]$. For this, it is sufficient to show that for every $A \in \mathcal{B}([0,T]) \otimes \mathcal{G}$ and $\alpha \in L^2_{\mathcal{G}}(\Omega;\mathbb{R}^m)$, the mapping $t \mapsto \mathbb{E}\left[ \mathbbm{1}_A(t,\omega) \alpha(\omega) \right]$ is Lebesgue measurable. By fixing $\alpha \in L^2_{\mathcal{G}}(\Omega;\mathbb{R}^m)$, this can be achieved by proving that the family $\{ A \in \mathcal{B}([0,T]) \otimes \mathcal{G} : \ t \mapsto \mathbb{E}\left[ \mathbbm{1}_A(t,\omega) \alpha(\omega) \right] \ \textnormal{is Lebesgue measurable} \}$ is a monotone class and then using standard monotone class arguments (the details are left to the reader). At this step, the Lebesgue differentiation theorem provides that for almost every $t \in [0,T]$, the following relations hold:
$$
\underset{\eta \rightarrow 0}{\lim} \ \frac{1}{\eta} \int^{t + \eta}_t \mathbb{E}\Big[ \| u(s) - u(t) \| \Big] \; \mathrm{d}s = 0 , \quad \underset{\eta \rightarrow 0}{\lim} \ \frac{1}{\eta} \int^{t + \eta}_t \mathbb{E}\Big[ \| u(s) - u(t) \|^2 \Big] \; \mathrm{d}s = 0 .
$$
Such a time $t \in [0,T]$ is called Lebesgue point for the control $u \in L^2_{\mathcal{F}}([0,T] \times \Omega;\mathbb{R}^m)$.

\subsubsection{Modified Needle-like Variations and End-point Mapping}

We may extend the concept of needle-like variations and of end-point mapping previously introduced in Section \ref{sec_needleLike} to the setting of free final time and of stochastic controls as follows.

Given an integer $j \in \mathbb{N}$, fix $j$ times $0 < t_1 < \dots < t_j < t_f$ which are Lebesgue points for $u$, and fix $j$ random variables $u_1,\dots,u_j$ such that $u_i \in L^2_{\mathcal{F}_{t_i}}(\Omega;U)$. For fixed scalars $0 \le \eta_i < t_{i+1} - t_i$, $i=1,\dots,j-1$, and $0 \le \eta_j < t_f - t_j$, the needle-like variation $\pi = \{ t_i , \eta_i , u_i \}_{i = 1,\dots,j}$ of the control $u$ is defined to be the admissible control $u_{\pi}(t) = u_i$ if $t \in [t_i,t_i+\eta_i]$ and $u_{\pi}(t) = u(t)$ otherwise. Denote by $\tilde x_v$ the solution related to an admissible control $v$ of the augmented system 
$$
\begin{cases}
    \mathrm{d}x(t) = b(t,v(t),x(t)) \; \mathrm{d}t + \sigma(t,x(t)) \; \mathrm{d}B_t , \quad x(0) = x^0 \medskip \\
    \mathrm{d}x^{n+1}(t) = f^0(t,v(t),x(t)) \; \mathrm{d}t , \hspace{14ex} x^{n+1}(0) = 0
\end{cases}
$$
and define the mapping $\tilde g : \mathbb{R}^{n+1} \rightarrow \mathbb{R}^{q+1} : \tilde x \mapsto (g(x),x^{n+1})$. For every fixed time $t \in (t_j,t_f]$, by denoting $\delta_t \triangleq \min \{ t_{i+1} - t_i , t - t_j , T - t : i = 1,\dots,j-1 \} > 0$, the end-point mapping at time $t$ is defined to be
\begin{align*}
    F^j_t : \ &\mathcal{C}^j_t \triangleq B^{j+1}_{\delta_t}(0) \cap ( \mathbb{R} \times \mathbb{R}^j_+ ) \rightarrow \mathbb{R}^{q+1} \\
    &(\delta,\eta_1,\dots,\eta_j) \mapsto \mathbb{E}\left[ \tilde g(\tilde x_{u_{\pi}}(t+\delta)) \right] - \mathbb{E}\left[ \tilde g(\tilde x_u(t)) \right] \nonumber
\end{align*}
where $B^{j+1}_{\rho}$ is the open ball in $\mathbb{R}^{j+1}$ of radius $\rho > 0$. Variations on the variable $\delta$ are necessary only if free-final-time problems are considered, in which case $(A'_2)$, and in particular the fact that $g$ is an affine function, play a crucial role for computations. Due to Lemma \ref{lemma_bound} (and to $(A'_2)$ in the case of free-final-time problems), it is not difficult to see that $F^j_t$ is Lipschitz (see also the argument developed to prove Lemma \ref{lemma_needleGeneralized} below). In addition, this mapping may be Gateaux differentiated at zero along admissible directions of the cone $\mathcal{C}^j_t$. For this, denote $\tilde b = (b^{\top},f^0)^{\top}$, $\tilde \sigma = (\sigma^{\top},0)^{\top}$ and let $z_{t_i,u_i}$ be the unique solution to \eqref{eq_rightSDE} with $\xi_{t_i} = \tilde b(t_i,u_i,x(t_i)) - \tilde b(t_i,u(t_i),x(t_i))$.
\begin{lmm}[Generalized stochastic needle-like variation formula] \label{lemma_needleGeneralized}
Let $(\delta,\eta_1,\dots,\eta_j) \in \mathcal{C}^j_t$ and assume $(A_2)$ holds (in particular, $\tilde g$ is an affine function when $\delta \neq 0$). If $t > t_j$ is a Lebesgue point for $u$, then it holds that
$$
\bigg\| \mathbb{E}\bigg[ \tilde g(\tilde x_{u_{\pi}}(t + \delta)) - \tilde g(\tilde x_u(t)) - \delta \frac{\partial \tilde g}{\partial \tilde x}(\tilde x_u(t)) \tilde b(t,u(t),x_u(t)) - \sum^j_{i=1} \eta_i \frac{\partial \tilde g}{\partial \tilde x}(\tilde x_u(t)) z_{t_i,u_i}(t) \bigg] \bigg\| = o\left( \delta + \sum^j_{i=1} \eta_i \right) .
$$
\end{lmm}
As for Lemma \ref{lemma_needle}, we provide an extensive proof of Lemma \ref{lemma_needleGeneralized} in the appendix.

\subsubsection{Variational Inequalities and Conclusion} \label{sec_variational}

Similarly to what we have argued in Section \ref{sec_variational}, the main step in the proof of the stochastic PMP with free final time and stochastic control goes by contradiction, leveraging Lemma \ref{lemma_needleGeneralized}. For this, from now on we assume that when free, the final time is a Lebesgue point for the optimal control. Otherwise, one may proceed by mimicking the argument developed in \cite[Section 7.3]{gamkrelidze2013}.

Assume that $(A_1)$ and $(A'_2)$ hold. For every $j \in \mathbb{N}$, define the linear mapping
$$
dF^j_{t_f}(\delta,\eta) = \delta \mathbb{E}\bigg[ \frac{\partial \tilde g}{\partial \tilde x}(\tilde x_u(t_f)) \tilde b(t_f,u(t_f),x_u(t_f)) \bigg] + \sum_{i=1}^j \eta_i \mathbb{E}\left[ \frac{\partial \tilde g}{\partial \tilde x}(\tilde x_u(t_f)) z_{t_i,u_i}(t_f) \right] ,
$$
which due to Lemma \ref{lemma_needleGeneralized}, satisfies
$$
\underset{\alpha > 0 , \alpha \rightarrow 0}{\lim} \ \frac{F^j_{t_f}\left( \alpha (\delta,\eta) \right)}{\alpha} = dF^j_{t_f}(\delta,\eta) ,
$$
for every $(\delta,\eta) \in \mathbb{R} \times \mathbb{R}^j_+$. Finally, consider the closed, convex cone of $\mathbb{R}^{q+1}$ given by
\begin{align*}
    K \triangleq \textnormal{Cl} \bigg( \textnormal{Cone} \ \bigg\{ &\mathbb{E}\left[ \frac{\partial \tilde g}{\partial \tilde x}(\tilde x_u(t_f)) z_{t_i,u_i}(t_f) \right] , \pm \mathbb{E}\left[ \frac{\partial \tilde g}{\partial \tilde x}(\tilde x_u(t_f)) \tilde b(t_f,u(t_f),x_u(t_f)) \right] \\ & \textnormal{for} \ u_i \in U \ \textnormal{and} \ t_i \in (0,t_f) \ \textnormal{is Lebesgue for} \ u \bigg\} \bigg) .
\end{align*}
If $K = \mathbb{R}^{q+1}$, it would hold $dF^j_{t_f}(\mathbb{R} \times \mathbb{R}^j_+) = K = \mathbb{R}^{q+1}$, and by \cite[Lemma 12.1]{agrachev2013}, one would find that the origin is an interior point of $F^j_{t_f}(\mathcal{C}^j_{t_f})$. This would imply that $(t_f,u,x)$ is not optimal for OCP, a contradiction.

The argument above (together with an application of the separation plane theorem) provides the existence of a non-zero vector denoted $\tilde{\mathfrak{p}} = (\mathfrak{p}^{\top},\mathfrak{p}^0) \in \mathbb{R}^{q+1}$ such that the following variational inequalities hold
{
\begin{equation} \label{eq_variationFinalTime}
\begin{cases}
    \displaystyle \tilde{\mathfrak{p}}^{\top} \mathbb{E}\left[ \frac{\partial \tilde g}{\partial \tilde x}(\tilde x_u(t_f)) \tilde b(t_f,u(t_f),x_u(t_f)) \right] = 0 \bigskip \\
    \displaystyle \tilde{\mathfrak{p}}^{\top} \mathbb{E}\left[ \frac{\partial \tilde g}{\partial \tilde x}(\tilde x_u(t_f)) z_{r,v}(t_f) \right] \le 0 , \ r \in [0,t_f] \ \textnormal{is Lebesgue for} \ u , \ v \in L^2_{\mathcal{F}_r}(\Omega;U) .
\end{cases}
\end{equation}}
In the case of deterministic controls, the random variables $v \in L^2_{\mathcal{F}_r}(\Omega;U)$ in \eqref{eq_variationFinalTime} are replaced by deterministic vectors $v \in U$. Moreover, when $t_f = T$, only negative variations on the final time are allowed. Hence in this case, the first equality of \eqref{eq_variationFinalTime} actually becomes a greater-or-equal-to-zero inequality (see also \cite[Chapter 7]{gamkrelidze2013}).

\vspace{5pt}

The rest of the proof remains unchanged, i.e., it follows the argument of Section \ref{sec_daulDerivation} verbatim, and the proof of Theorem \ref{theo_mainGeneralized} may be readily inferred by the computations of Section \ref{sec_proofConv}. 

\section{An Example Numerical Scheme} \label{sec_numerical}

Although the procedure detailed previously provides methodological steps to tackle OCP through successive linearizations, numerically solving LCPs that depend on stochastic coefficients still remains a challenge. In this last section, under appropriate assumptions we propose an approximate, though very practical, numerical scheme to effectively solve each subproblem LOCP$^{\Delta}_k$. We stress the fact that our main goal is not the development of an ultimate algorithm, but rather consists of demonstrating how one may leverage the theoretical insights provided by Theorem \ref{theo_main} to design efficient strategies to practically solve OCP.

\subsection{A Simplified Context}\label{sec:num:simplified_context} The proposed approach relies on a specific shape of the cost and the dynamics of OCP. Specifically, we consider OCPs with deterministic admissible controls $u \in \mathcal{U} = L^2([0,T];U)$ (over a fixed time horizon $[0,T]$) only and whose cost functions $f^0$ are such that $H = 0$. Moreover, we assume the state variable is given by two components $(x,z) \in \mathbb{R}^{n_x+n_z}$ for $n_x, n_z \in \mathbb{N}$, satisfying the following system of forward stochastic differential equations ($b^x$ and $b^z$ are accordingly defined as in \eqref{eq_SDE})
\begin{equation} \label{eq_SDENumerics}
    \begin{cases}
        \mathrm{d}x(t) = b^x(t,u(t),x(t),z(t)) \; \mathrm{d}t + \sigma(t,z(t)) \; \mathrm{d}B_t , \quad x(0) = x^0 \medskip \\
        \mathrm{d}z(t) = b^z(t,u(t),z(t)) \; \mathrm{d}t , \quad z(0) = z^0 \in \mathbb{R}^{n_z} .
    \end{cases}
\end{equation}
In particular, any $\mathcal{F}$-adapted process $(x,z)$ solution to \eqref{eq_SDENumerics} with continuous sample paths for a given control $u \in \mathcal{U}$ is such that $z$ is deterministic. For the sake of clarity and to avoid cumbersome notation, from now on we assume $g = 0$ and that $b^x$ does not explicitly depend on $z$. This is clearly done without loss of generality.

\subsection{The Proposed Approach} With the assumptions adopted previously, we see that the diffusion in the dynamics of OCP is now forced to be deterministic. This fact is at the root of our method, which mimics the procedure proposed in \cite{berret2020}. Specifically, we transcribe every stochastic subproblem LOCP$^{\Delta}_k$ into a deterministic and convex optimal control problem, whereby the variables are the mean and the covariance of the solution $x_k$ to LOCP$^{\Delta}_k$. The main advantage in doing so is that deterministic reformulations of the subproblems LOCP$^{\Delta}_k$ can be efficiently solved via off-the-shelf convex solvers. Unlike \cite{berret2020}, here we rely on some \textit{upstream information} for the design of this numerical scheme which is entailed by Theorem \ref{theo_main} as follows.

By recalling the notation introduced in the previous sections, we denote $\mu(t) \triangleq \mathbb{E}[x(t)]$, $\Sigma(t) \triangleq \mathbb{E}[(x(t) - \mu(t)) (x(t) - \mu(t))^{\top}]$, and for $k \in \mathbb{N}$, $\mu_k(t) \triangleq \mathbb{E}[x_k(t)]$ and $\Sigma_k(t) \triangleq \mathbb{E}[(x_k(t) - \mu_k(t)) (x_k(t) - \mu_k(t))^{\top}]$. Heuristically, assuming that solutions to LOCP$^{\Delta}_k$ have \textit{small variance}, i.e.,  $\| \textnormal{tr} \ \Sigma_k \|_{L^2} \ll 1$, one can compute the linearization of $b^x$, $b^z$, $\sigma$, and $L$ at $\mu_k$ rather than at $x_k$. In doing so, for the cost of LOCP$^{\Delta}_{k+1}$ we obtain the following approximation (the notation goes accordingly as in \eqref{eq_LCost})
\begin{equation} \label{eq_LCostNumerics}
    \mathbb{E}\left[ \int^T_0 f^0_{k+1}(s,u(s),x(s),z(s)) \; \mathrm{d}s \right] \approx \int^T_0 f^0_{u_k,\mu_k}(s,u(s),\mu(s),z(s)) \; \mathrm{d}s ,
\end{equation}
whereas for the dynamics of LOCP$^{\Delta}_{k+1}$ (the notation goes accordingly as in \eqref{eq_LSDE})
\begin{equation} \label{eq_LSDENumerics}
    \begin{cases}
        \mathrm{d}x(t) \approx b^x_{u_k,\mu_k}(t,u(t),x(t)) \; \mathrm{d}t + \sigma_{z_k}(t,z(t)) \; \mathrm{d}B_t \medskip \\
        \hspace{6ex} = \Big( \mathcal{A}_{k+1}(t) x(t) + \mathcal{B}_{k+1}(t,u(t)) \Big) \; \mathrm{d}t + \mathcal{C}_{k+1}(t,z(t)) \; \mathrm{d}B_t \medskip \\
        \mathrm{d}z(t) \approx b^z_{u_k,z_k}(t,u(t),z(t)) \; \mathrm{d}t = \Big( \mathcal{D}_{k+1}(t) z(t) + \mathcal{E}_{k+1}(t,u(t)) \Big) \; \mathrm{d}t
    \end{cases}
\end{equation}
where $\mathcal{A}_{k+1}(t) \in \mathbb{R}^{n_x \times n_x}$, $\mathcal{B}_{k+1}(t,u)  \in \mathbb{R}^{n_x}$, $\mathcal{C}_{k+1}(t,z) \in \mathbb{R}^{n_x}$, $\mathcal{D}_{k+1}(t) \in \mathbb{R}^{n_z \times n_z}$, and $\mathcal{E}_{k+1}(t,u) \in \mathbb{R}^{n_z}$ are deterministic and with $\mathcal{C}_{k+1}(t,z)$ affine in $z$. Accordingly, by introducing $\mu^{e_k}(t) \triangleq \mu(t) - \mu_k(t)$ and $\Sigma^{e_k}(t) \triangleq \mathbb{E}[(x(t) - x_k(t) - \mu(t) + \mu_k(t)) (x(t) - x_k(t) - \mu(t) + \mu_k(t))^{\top}]$, slightly tighter trust-region constraints are
\begin{equation} \label{eq_trustApprox}
    \int^T_0 \textnormal{tr} \ \Sigma^{e_k}(t) \; \mathrm{d}t + \int^T_0 \| \mu^{e_k}(t) \|^2 \; \mathrm{d}t \le \Delta_{k+1} , \quad | t^{k+1}_f - t^k_f | \le \Delta_{k+1} .
\end{equation}
At this point, as all coefficients are deterministic, solutions to \eqref{eq_LSDENumerics} are Gaussian processes whose dynamics take the form (see, e.g., \cite{maybeck1982})
\begin{equation} \label{eq_LSDEmean}
    \begin{cases}
        \dot{\mu}(t) = \mathcal{A}_{k+1}(t) \mu(t) + \mathcal{B}_{k+1}(t,u(t)) , \quad \dot{z}(t) = \mathcal{D}_{k+1}(t) z(t) + \mathcal{E}_{k+1}(t,u(t)) \medskip \\
        \dot{\Sigma}(t) = \mathcal{A}_{k+1}(t) \Sigma(t) + \Sigma(t) \mathcal{A}_{k+1}(t)^{\top} + \mathcal{C}_{k+1}(t,z(t)) \mathcal{C}_{k+1}(t,z(t))^{\top} .
    \end{cases}
\end{equation}
The system above is not linear because of $\mathcal{C}_{k+1}(t,z(t)) \mathcal{C}_{k+1}(t,z(t))^{\top}$. Nevertheless, we may call upon the convergences of Theorem \ref{theo_main} 
, and in particular the convergence of the sequence $(z_k)_{k \in \mathbb{N}}$ to some deterministic curve $z$, to replace \eqref{eq_LSDEmean} with
\begin{equation} \label{eq_LSDEmeanFinal}
    \begin{cases}
        \dot{\mu}(t) = \mathcal{A}_{k+1}(t) \mu(t) + \mathcal{B}_{k+1}(t,u(t)) , \quad \dot{z}(t) = \mathcal{D}_{k+1}(t) z(t) + \mathcal{E}_{k+1}(t,u(t)) \medskip \\
        \dot{\Sigma}(t) = \mathcal{A}_{k+1}(t) \Sigma(t) + \Sigma(t) \mathcal{A}_{k+1}(t)^{\top} + \mathcal{C}_{k+1}(t,z_k(t)) \mathcal{C}_{k+1}(t,z(t))^{\top} .
    \end{cases}
\end{equation}
In conclusion, we may heuristically replace every LOCP$^{\Delta}_{k+1}$ with the deterministic convex optimal control problem  whose dynamics are \eqref{eq_LSDEmeanFinal}, whose trust-region constraints are \eqref{eq_trustApprox}, whose variables are $\mu$ and $\Sigma$ (and additionally $\mu^{e_k}$ and $\Sigma_{e_k}$), and whose cost consists of replacing \eqref{eq_LCostNumerics} with
$$
\int^T_0 \Big( f^0_{u_k,\mu_k}(s,u(s),\mu(s),z(s)) + \textnormal{tr} \ \Sigma(t) \Big) \; \mathrm{d}s
$$
to force solutions to have small variances, which in turn justifies the whole approach.

\subsection{Uncertain Car Trajectory Planning Problem} \label{sec:results}

We now provide numerical experiments. We consider the problem of planning the trajectory of a car whose state and control inputs are $x=(r_x,r_y,\theta, v,\omega)$ and $u=(a_v,a_{\omega})$ and whose dynamics are
\vspace{-2mm}
\begin{equation}
b(t,x,u)=
(v\cos(\theta),
v\sin(\theta),
\omega,
a_v,
a_{\omega}),
\
\sigma(t,x) = \textrm{diag}([
\alpha^2\omega v, \alpha^2\omega v, \beta^2\omega v, 0, 0]),
\nonumber
\end{equation}
where $\alpha^2=0.1$ and $\beta^2=0.01$ quantify the effect of slip. The evolution of $(v,\omega)$ is deterministic: given actuator commands, the change in velocity is known exactly, but uncertainty in positional variables persists. 
This model is motivated by driving applications where the vehicle may drift during tight turns and suffer from important positional uncertainty, but limits its acceleration to avoid wheel slip. 
By defining $z=(v,\omega)$, this problem setting matches the one presented in the previous section.
We consider minimizing control effort $G(u(s))=\|u(s)\|_2^2$.  
The initial state $x^0$ is fixed, and the final state constraint consists of reaching the goal $\smash{x^f}$ in expectation such that $\smash{g(x(t_f))=x(t_f)-x^f}$. 
Further, we 
consider cylindrical obstacles of radius $\delta_o>0$ centered at $r_o\in\mathbb{R}^2$ %
for which we define the 
potential function $c_o:\mathbb{R}^n\rightarrow\mathbb{R}$ as $c_o(x)=\| r - r_o \|^2 - (\delta_o+\varepsilon)^2$ if $\| r - r_o \| < (\delta_o+\varepsilon)$ and $0$ otherwise, 
where $r=[r_x,r_y]\subset x$ and $\varepsilon=0.1$ is an additional clearance. We penalize obstacle avoidance constraint violation directly within the cost, defining $L(s,x(s),u(s))=L_0(s,x(s))=\lambda c_o(x(s))$, with $\lambda = 500$ (note that $H=0$ in this setting, following Section \ref{sec:num:simplified_context}).

\subsection{Results}
\begin{figure}[!htb]
    \begin{minipage}{.5\linewidth}
        \centering
        \includegraphics[width=0.9\linewidth]{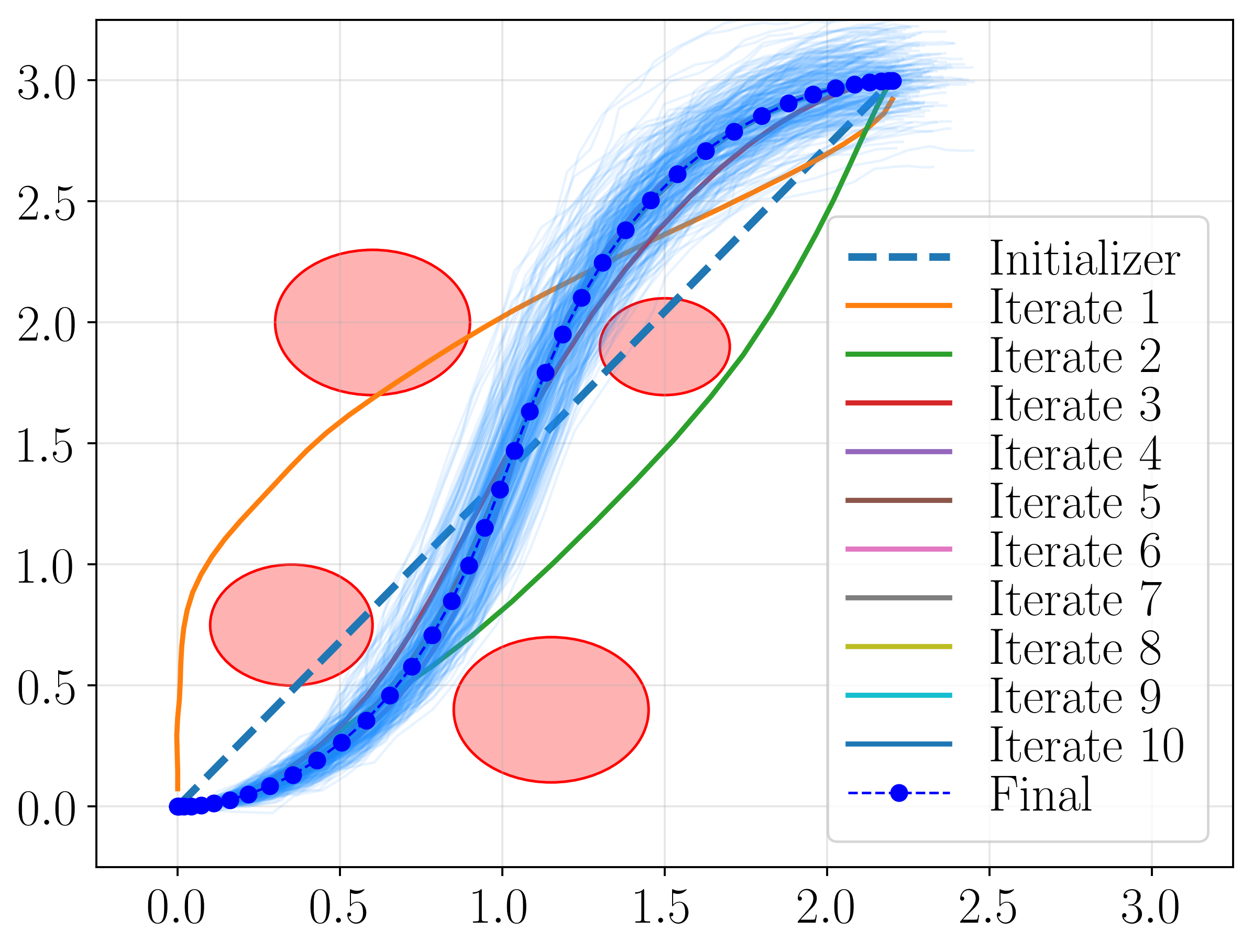}
    \end{minipage}%
    \begin{minipage}{.495\linewidth}
        \centering
        \includegraphics[width=0.45\linewidth]{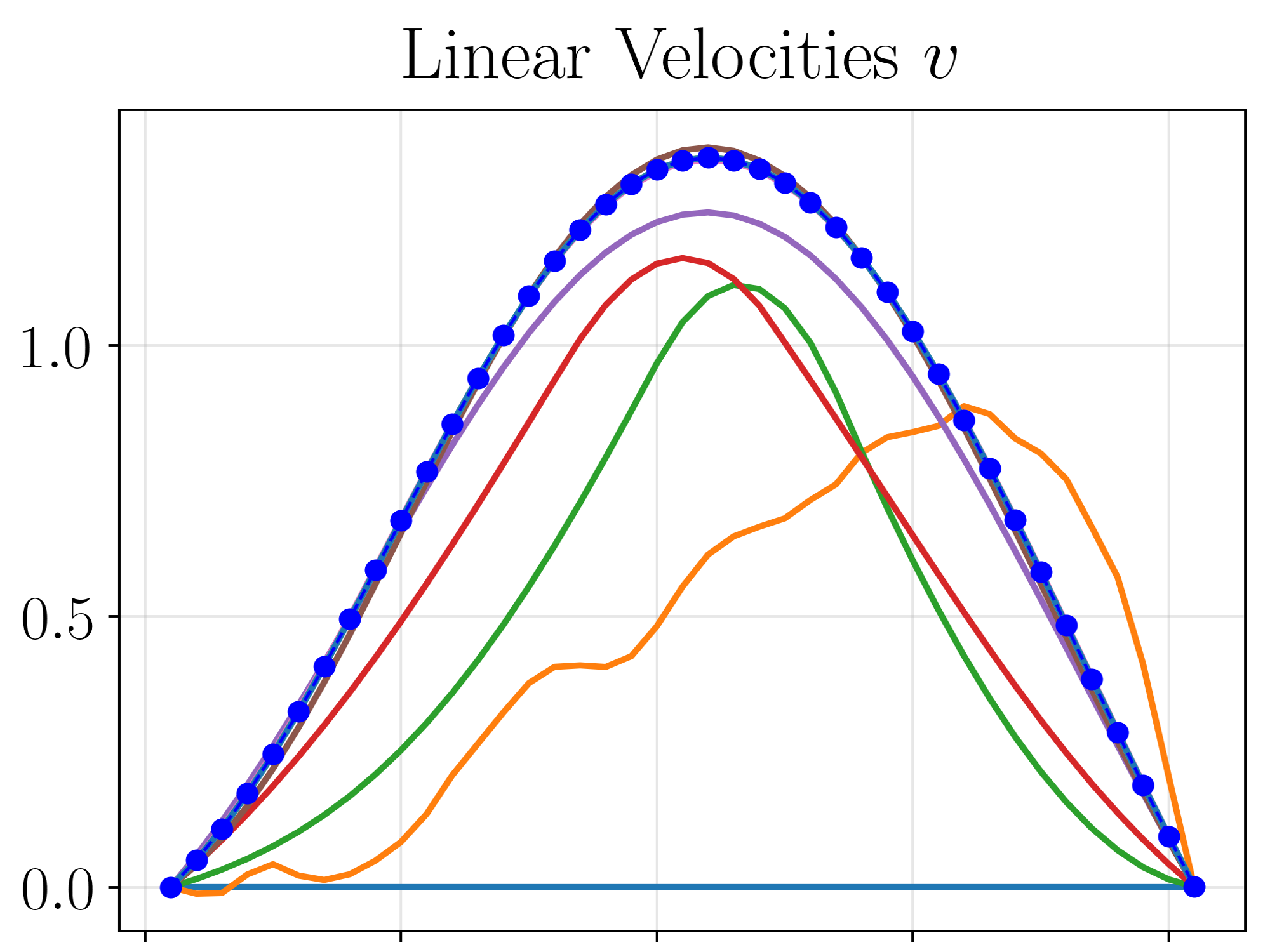}
        \includegraphics[width=0.45\linewidth]{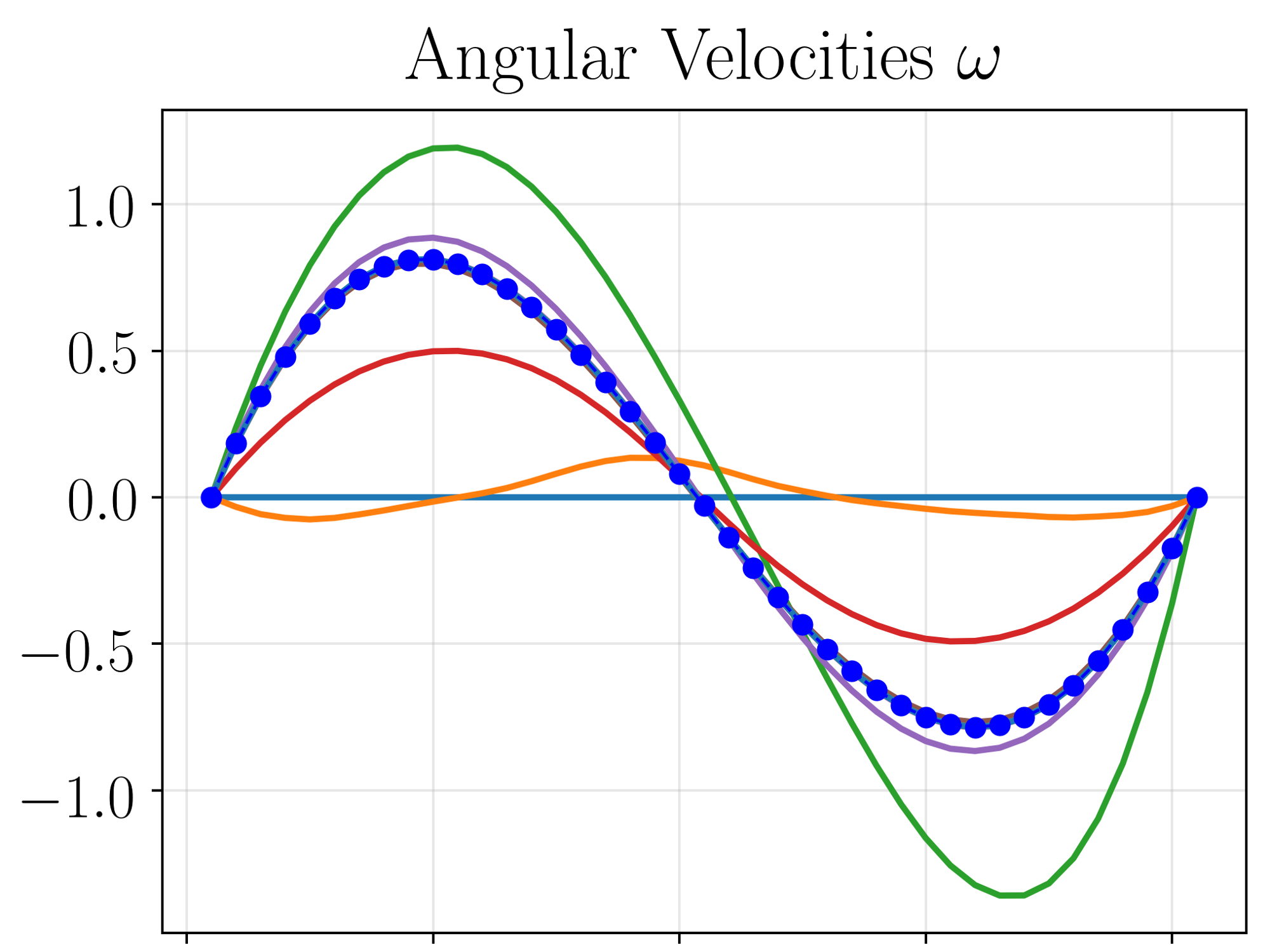}
        \\
        \includegraphics[width=0.45\linewidth]{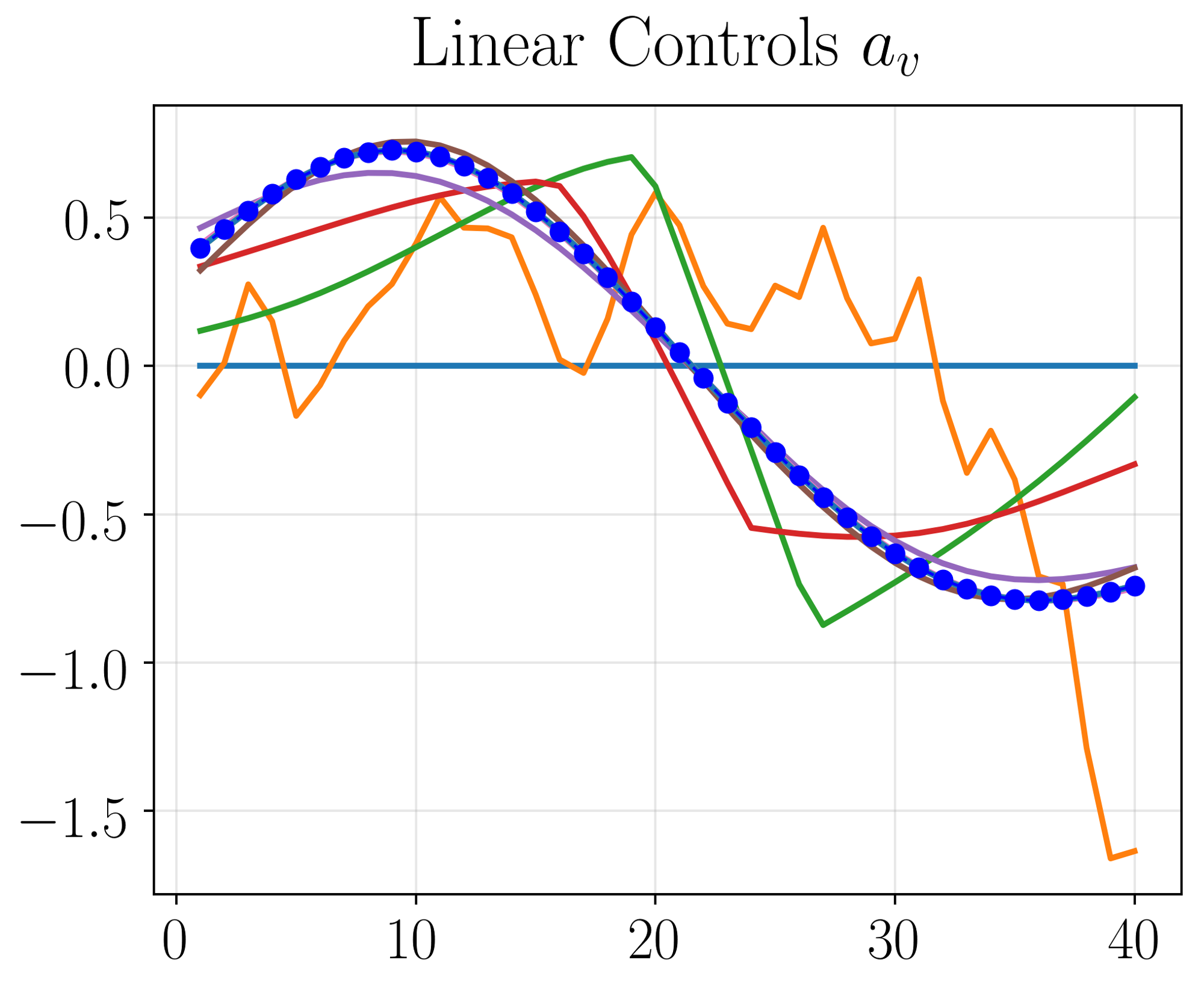}
        \includegraphics[width=0.45\linewidth]{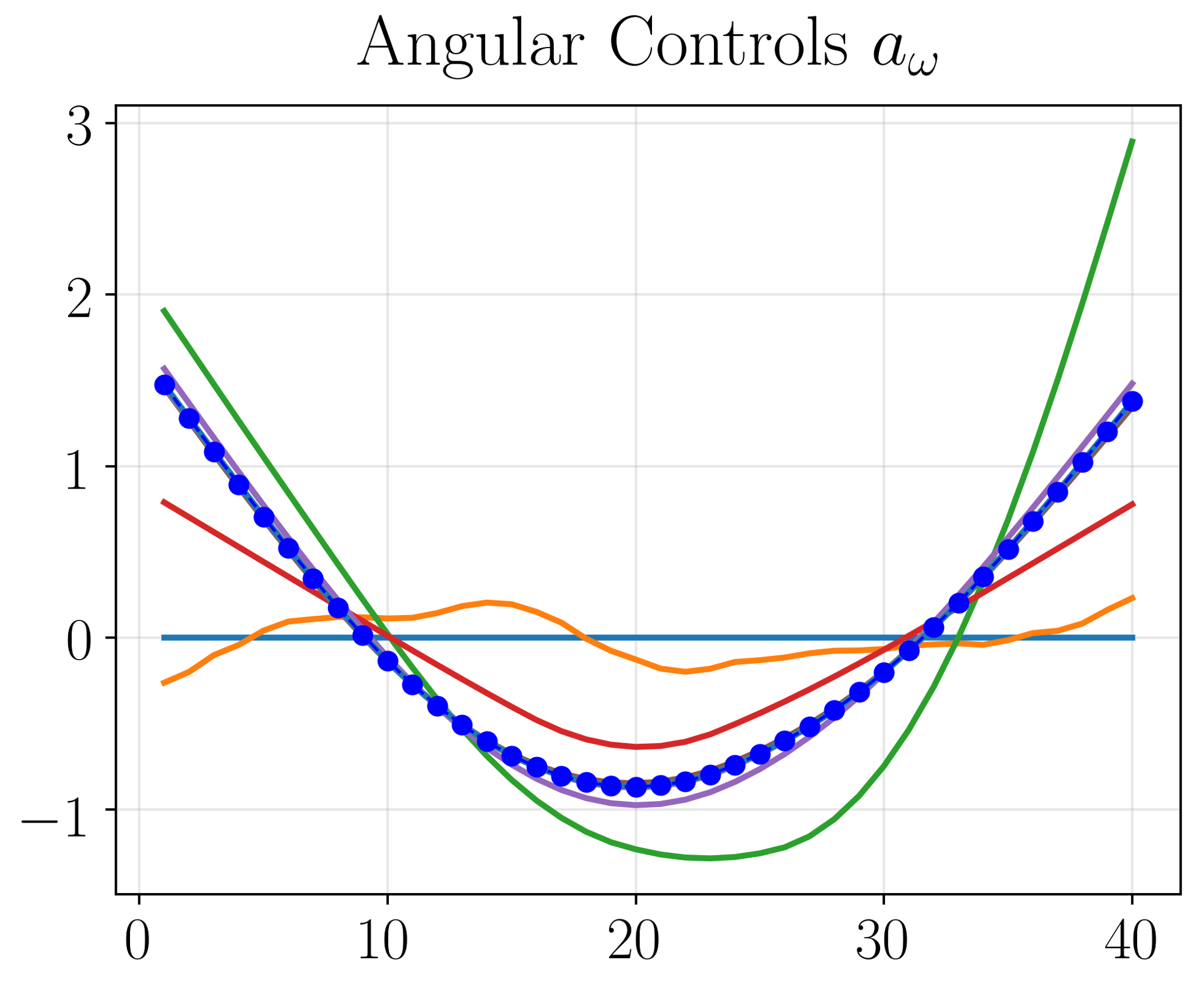}
    \end{minipage}%
    \caption{Left: solution of each SCP iteration and $100$ sample paths of the resulting final trajectory. 
    Right: velocities and control inputs at each SCP iteration.}
    \label{fig:results}
\end{figure}

We set $x(0)=[0,\dots,0]$, 
$x^f=[2.2,3,0,0,0]$, and 
$t_f=5~\textrm{s}$, and we set up four obstacles. 
We discretize \eqref{eq_LSDEmeanFinal} using a forward Euler scheme with $N=41$ discretization nodes, set $\Delta_0=100$, set $\Delta_{k+1}=0.99\Delta_k$ at each SCP iteration, %
and use IPOPT to solve each convexified problem.  
We check convergence of SCP by verifying that $\int^{t_f}_0 \| u_{k+1}-u_k \|^2(s) + \| u_{k}-u_{k-1} \|^2(s) \; \mathrm{d}s \leq 10^{-3}$. 
Although our method is initialized with an unfeasible straight-line initial trajectory, SCP converges in $10$ iterations, with a trajectory avoiding the obstacles in expectation. Indeed, after evaluating $10^4$ sample paths of the system, only $6.14\%$ of the trajectories intersect with obstacles. 
\firstRev{We also verify that $p^0_k\neq 0$ at every SCP iteration, and at convergence, so that $(\mathfrak{p},p^0) \neq 0$. Thus, Theorems \ref{theo_main} and \ref{theo_main2} apply and the solution generated by SCP is an extremal for OCP. We visualize our results in Figure \ref{fig:results} and release our implementation at \url{https://github.com/StanfordASL/stochasticSCP}}. 

\section{Conclusion and Perspectives} \label{sec_conclusion}

In this paper we introduced and analyzed convergence properties for sequential convex programming for non-linear stochastic optimal control, from which we derived a practical numerical framework to solve non-linear stochastic optimal control problems.

Future work may consider  
extending this analysis to tackle more general problem formulations, e.g., risk measures as costs and state (chance) constraints. In this context, some preliminary results using SCP only exist for discrete time problem formulations \cite{LewBonalli2020}. However, tackling continuous time formulations will require more sophisticated necessary conditions for optimality. 
In addition, we plan to better investigate the setting of free final time and stochastic admissible controls, thus devising sharper assumptions and improving the convergence result. Finally, we plan to further leverage our theoretical insights to design new and more efficient numerical schemes for non-linear stochastic optimal control. %

\bibliographystyle{unsrt}
\bibliography{ASL_papers,references}

\newcommand{\noopsort}[1]{} \newcommand{\printfirst}[2]{#1}
  \newcommand{\singleletter}[1]{#1} \newcommand{\switchargs}[2]{#2#1}
\begin{thebibliography}{10}

\bibitem{potter1965}
J.E. Potter.
\newblock A matrix equation arising in statistical filter theory.
\newblock {\em Rep. RE-9, Experimental Astronomy Laboratory, Massachusetts
  Institute of Technology}, 1965.

\bibitem{bismut1976}
J.-M. Bismut.
\newblock Linear quadratic optimal stochastic control with random coefficients.
\newblock {\em SIAM Journal on Control and Optimization}, 14(3):419--444, 1976.

\bibitem{peng1992}
S.~Peng.
\newblock Stochastic {Hamilton--Jacobi--Bellman} equations.
\newblock {\em SIAM Journal on Control and Optimization}, 30(2):284--304, 1992.

\bibitem{tang2003}
S.~Tang.
\newblock General linear quadratic optimal stochastic control problems with
  random coefficients: linear stochastic {Hamilton} systems and backward
  stochastic {Riccati} equations.
\newblock {\em SIAM journal on control and optimization}, 42(1):53--75, 2003.

\bibitem{kleindorfer1973}
P.~Kleindorfer and K.~Glover.
\newblock Linear convex stochastic optimal control with applications in
  production planning.
\newblock {\em IEEE Transactions on Automatic Control}, 18(1):56--59, 1973.

\bibitem{rockafellar1990}
R.~T. Rockafellar and R.~J.~B. Wets.
\newblock Generalized linear-quadratic problems of deterministic and stochastic
  optimal control in discrete time.
\newblock {\em SIAM Journal on control and optimization}, 28(4):810--822, 1990.

\bibitem{kuhn2011}
D.~Kuhn, W.~Wiesemann, and A.~Georghiou.
\newblock Primal and dual linear decision rules in stochastic and robust
  optimization.
\newblock {\em Mathematical Programming}, 130(1):177--209, 2011.

\bibitem{bes1989}
C.~Bes and S.~Sethi.
\newblock Solution of a class of stochastic linear-convex control problems
  using deterministic equivalents.
\newblock {\em Journal of optimization theory and applications}, 62(1):17--27,
  1989.

\bibitem{berret2020}
B.~Berret and F.~Jean.
\newblock Efficient computation of optimal open-loop controls for stochastic
  systems.
\newblock {\em Automatica}, 115(108874), 2020.

\bibitem{rami2000}
M.~A. Rami and X.~Y. Zhou.
\newblock Linear matrix inequalities, riccati equations, and indefinite
  stochastic linear quadratic controls.
\newblock {\em IEEE Transactions on Automatic Control}, 45(6):1131--1143, 2000.

\bibitem{yao2001}
D.~D. Yao, S.~Zhang, and X.~Y. Zhou.
\newblock Stochastic linear-quadratic control via semidefinite programming.
\newblock {\em SIAM Journal on Control and Optimization}, 40(3):801--823, 2001.

\bibitem{Bertsimas2007}
D.~Bertsimas and D.~B. Brown.
\newblock Constrained stochastic lqc: a tractable approach.
\newblock {\em IEEE Transactions on automatic control}, 52(10):1826--1841,
  2007.

\bibitem{damm2017}
T.~Damm, H.~Mena, and T.~Stillfjord.
\newblock Numerical solution of the finite horizon stochastic linear quadratic
  control problem.
\newblock {\em Numerical Linear Algebra with Applications}, 24(4), 2017.

\bibitem{levajkovic2018}
T.~Levajkovi{\'c}, H.~Mena, and L.-M. Pfurtscheller.
\newblock Solving stochastic {LQR} problems by polynomial chaos.
\newblock {\em IEEE control systems letters}, 2(4):641--646, 2018.

\bibitem{bellman1957}
R.~Bellman.
\newblock {\em Dynamic Programming}.
\newblock Princeton Univ. Press, Princeton, New Jersey, 1957.

\bibitem{lions1983}
P.-L. Lions.
\newblock Optimal control of diffusion processes and {Hamilton-Jacobi-Bellman}
  equations, part i.
\newblock {\em Comm. Partial Differential Equations}, 8:1101--1174, 1983.

\bibitem{pontryagin2018}
L.S. Pontryagin.
\newblock {\em Mathematical theory of optimal processes}.
\newblock Routledge, 2018.

\bibitem{kushner1964}
H.~J. Kushner and F.~C. Schweppe.
\newblock A maximum principle for stochastic control systems.
\newblock {\em Journal of Mathematical Analysis and Applications},
  8(2):287--302, 1964.

\bibitem{Peng1990}
S.~Peng.
\newblock A general stochastic maximum principle for optimal control problems.
\newblock {\em SIAM Journal on control and optimization}, 28(4):966--979, 1990.

\bibitem{Yong1999}
J.~Yong and X.~Y. Zhou.
\newblock {\em Stochastic controls: {Hamiltonian} systems and {HJB} equations},
  volume~43.
\newblock Springer Science \& Business Media, 1999.

\bibitem{trelat2012}
E.~Tr{\'e}lat.
\newblock Optimal control and applications to aerospace: some results and
  challenges.
\newblock {\em Journal of Optimization Theory and Applications},
  154(3):713--758, 2012.

\bibitem{Bonalli2017Bis}
R.~Bonalli, B.~H{\'e}riss{\'e}, and E.~Tr{\'e}lat.
\newblock {Analytical Initialization of a Continuation-Based Indirect Method
  for Optimal Control of Endo-Atmospheric Launch Vehicle Systems}.
\newblock In {\em IFAC World Congress}, 2017.

\bibitem{bonalli2019_TAC}
R.~Bonalli, B.~H{\'e}riss{\'e}, and E.~Tr{\'e}lat.
\newblock Optimal control of endo-atmospheric launch vehicle systems: Geometric
  and computational issues.
\newblock {\em IEEE Transactions on Automatic Control}, 65(6):2418--2433, 2019.

\bibitem{shapiro2005}
A.~Shapiro and A.~Nemirovski.
\newblock On complexity of stochastic programming problems.
\newblock In {\em Continuous optimization}, pages 111--146. Springer, 2005.

\bibitem{gobet2016}
E.~Gobet.
\newblock {\em Monte-Carlo methods and stochastic processes: from linear to
  non-linear}.
\newblock CRC Press, 2016.

\bibitem{kushner1990}
H.~J. Kushner.
\newblock Numerical methods for stochastic control problems in continuous time.
\newblock {\em SIAM Journal on Control and Optimization}, 28(5):999--1048,
  1990.

\bibitem{kushner1991}
H.~J. Kushner and L.~F. Martins.
\newblock Numerical methods for stochastic singular control problems.
\newblock {\em SIAM journal on control and optimization}, 29(6):1443--1475,
  1991.

\bibitem{annunziato2013}
M.~Annunziato and A.~Borz{\`\i}.
\newblock A {Fokker}--{Planck} control framework for multidimensional
  stochastic processes.
\newblock {\em Journal of Computational and Applied Mathematics},
  237(1):487--507, 2013.

\bibitem{LeGall2016}
J.-F. Le~Gall.
\newblock {\em Brownian motion, martingales, and stochastic calculus}, volume
  274.
\newblock Springer, 2016.

\bibitem{Carmona2016}
R.~Carmona.
\newblock {\em Lectures on BSDEs, stochastic control, and stochastic
  differential games with financial applications}, volume~1.
\newblock SIAM, 2016.

\bibitem{Kazuhide2018}
O.~Kazuhide, M.~Goldshtein, and P.~Tsiotras.
\newblock Optimal covariance control for stochastic systems under chance
  constraints.
\newblock {\em {IEEE Control Systems Letters}}, 2(2):266--271, 2018.

\bibitem{LewBonalli2020}
T.~Lew, R.~Bonalli, and M.~Pavone.
\newblock Chance-constrained sequential convex programming for robust
  trajectory optimization.
\newblock In {\em {European Control Conference}}, 2020.

\bibitem{wang2017}
Y.~Wang, D.~Yang, J.~Yong, and Z.~Yu.
\newblock Exact controllability of linear stochastic differential equations and
  related problems.
\newblock {\em American Institute of Mathematical Sciences}, 7(2):305--345,
  2017.

\bibitem{Dinh2010}
Q.~T. Dinh and M.~Diehl.
\newblock {Local Convergence of Sequential Convex Programming for Nonconvex
  Optimization}.
\newblock In {\em {Recent Advances in Optimization and its Applications in
  Engineering}}, pages 93--102. Springer, 2010.

\bibitem{Diehl2019}
M.~Diehl and F.~Messerer.
\newblock {Local Convergence of Generalized Gauss-Newton and Sequential Convex
  Programming}.
\newblock In {\em Conference on Decision and Control}, 2019.

\bibitem{Mao2017}
Y.~Mao, D.~Dueri, M.~Szmuk, and B.~A\c{c}ikme\c{s}e.
\newblock {Successive Convexification of Non-Convex Optimal Control Problems
  with State Constraints}.
\newblock {\em {IFAC-PapersOnLine}}, 50(1):4063--4069, 2017.

\bibitem{palacios1982}
F.~Palacios-Gomez, L.~Lasdon, and M.~Engquist.
\newblock Nonlinear optimization by successive linear programming.
\newblock {\em Management science}, 28(10):1106--1120, 1982.

\bibitem{boggs1995}
P.~T. Boggs and J.~W. Tolle.
\newblock Sequential quadratic programming.
\newblock {\em Acta numerica}, 4(1):1--51, 1995.

\bibitem{gamkrelidze2013}
R.~Gamkrelidze.
\newblock {\em Principles of optimal control theory}, volume~7.
\newblock Springer Science \& Business Media, 2013.

\bibitem{agrachev2013}
A.A. Agrachev and Y.~Sachkov.
\newblock {\em Control theory from the geometric viewpoint}, volume~87.
\newblock Springer Science \& Business Media, 2013.

\bibitem{Frankowska2018}
H.~Frankowska, H.~Zhang, and X.~Zhang.
\newblock {Stochastic Optimal Control Problems with Control and Initial-Final
  States Constraints}.
\newblock {\em SIAM Journal on Control and Optimization}, 56:1823--1855, 2018.

\bibitem{Nocedal1999}
J.~Nocedal and S.~Wright.
\newblock {\em {Numerical Optimization}}.
\newblock Springer, 1999.

\bibitem{Lu2013}
Z.~Lu.
\newblock Sequential convex programming methods for a class of structured
  nonlinear programming.
\newblock Technical report, 2013.

\bibitem{Haberkorn2011}
T.~Haberkorn and E.~Tr\'elat.
\newblock Convergence results for smooth regularizations of hybrid nonlinear
  optimal control problems.
\newblock {\em SIAM Journal on Control and Optimization}, 49:1498--1522, 2011.

\bibitem{Bonalli2017}
R.~Bonalli, B.~H{\'e}riss{\'e}, and E.~Tr{\'e}lat.
\newblock {Solving Optimal Control Problems for Delayed Control-Affine Systems
  with Quadratic Cost by Numerical Continuation}.
\newblock In {\em American Control Conference}, 2017.

\bibitem{Bonalli2018}
R.~Bonalli.
\newblock {\em Optimal control of aerospace systems with control-state
  constraints and delays}.
\newblock PhD thesis, {Sorbonne Universit{\'e}}, 2018.

\bibitem{bonalli2019}
R.~Bonalli, B.~H{\'e}riss{\'e}, and E.~Tr\'elat.
\newblock Continuity of pontryagin extremals with respect to delays in
  nonlinear optimal control.
\newblock {\em SIAM Journal on Control and Optimization}, 57(2):1440--1466,
  2019.

\bibitem{ShvartsmanVinter2006}
I.~A. Shvartsman and R.~B. Vinter.
\newblock Regularity properties of optimal controls for problems with
  time-varying state and control constraints.
\newblock {\em Nonlinear Analysis: Theory, Methods \& Applications},
  65:448--474, 2006.

\bibitem{ChitourJeanEtAl2008}
Y.~Chitour, F.~Jean, and E.~{Tr\'elat}.
\newblock Singular trajectories of control-affine systems.
\newblock {\em SIAM Journal on Control and Optimization}, 47:1078--1095, 2008.

\bibitem{bryson1975}
A.E. Bryson.
\newblock {\em Applied optimal control: optimization, estimation and control}.
\newblock CRC Press, 1975.

\bibitem{betts1998}
J.T. Betts.
\newblock Survey of numerical methods for trajectory optimization.
\newblock {\em Journal of guidance, control, and dynamics}, 21(2):193--207,
  1998.

\bibitem{trelat2000}
Emmanuel Tr{\'e}lat.
\newblock Some properties of the value function and its level sets for affine
  control systems with quadratic cost.
\newblock {\em Journal of Dynamical and Control Systems}, 6(4):511--541, 2000.

\bibitem{bourdin2013}
L.~Bourdin and E.~Tr\'elat.
\newblock Pontryagin maximum principle for finite dimensional nonlinear optimal
  control problem on time scales.
\newblock {\em SIAM Journal on Control and Optimization}, 51(5):3781--3813,
  2013.

\bibitem{maybeck1982}
P.~S. Maybeck.
\newblock {\em Stochastic models, estimation, and control}.
\newblock Academic press, 1982.

\end{thebibliography}

\section*{Appendix}

\subsection{Proof of Lemma \ref{lemma_bound}} \label{sec_app1}

\begin{proof}[Proof of Lemma \ref{lemma_bound}]
For the sake of clarity in the notation, we only consider the case for which $d = 1$, the case with multivariate Brownian motion being similar.

Let us start with the first inequality. For this, by denoting $x = x_{u,v,y}$, for every $t \in [0,T]$ we compute (below, $C \ge 0$ is an appropriate constant)
{\small
\begin{align*}
    &\mathbb{E}\left[ \underset{s \in [0,t]}{\sup} \ \| x(s) \|^{\ell} \right] \le C \mathbb{E}\left[ \| x^0 \|^{\ell} \right] + C \mathbb{E}\left[ \underset{r \in [0,t]}{\sup} \ \left\| \int^r_0 \sigma(s,y(s)) \; \mathrm{d}B_s \right\|^{\ell} \right] \\
    & + C \mathbb{E}\left[ \underset{r \in [0,t]}{\sup} \ \left\| \int^r_0 \left( b_0(s,y(s)) + \sum^m_{i=1} u^i(s) b_i(s,y(s)) \right) \; \mathrm{d}s \right\|^p \right] \\
    & + C \mathbb{E}\left[ \underset{r \in [0,t]}{\sup} \ \left\| \int^r_0 \left( \frac{\partial b_0}{\partial x}(s,y(s)) + \sum^m_{i=1} v^i(s) \frac{\partial b_i}{\partial x}(s,y(s)) \right) (x(s) - y(s)) \; \mathrm{d}s \right\|^{\ell} \right] \\
    & + C \mathbb{E}\left[ \underset{r \in [0,t]}{\sup} \ \left\| \int^r_0 \frac{\partial \sigma}{\partial x}(s,y(s)) (x(s) - y(s)) \; \mathrm{d}B_s \right\|^{\ell} \right] .
\end{align*}}
For the last term, by denoting $S_{y,\sigma} \triangleq \bigg\{ (s,\omega) \in [0,T] \times \Omega : \ (s,y(s,\omega)) \in \textnormal{supp} \ \sigma \bigg\}$, Burkholder--Davis--Gundy, H\"older, and Young inequalities give
\begin{align*}
    &\mathbb{E}\left[ \underset{r \in [0,t]}{\sup} \ \left\| \int^r_0 \frac{\partial \sigma}{\partial x}(s,y(s)) (x(s) - y(s)) \; \mathrm{d}B_s \right\|^p \right] \\
    & \le C \left( \int^t_0 \mathbb{E}\left[ \underset{s \in [0,r]}{\sup} \ \left\| x(r) \right\|^p \right] \; \mathrm{d}s + \int_{S_{y,\sigma}} \left\| \frac{\partial \sigma}{\partial x}(s,y(s)) \right\|^p \left\| y(s) \right\|^p \; \mathrm{d}(s \times P) \right) \\
    & \le C \left( 1 + \int^t_0 \mathbb{E}\left[ \underset{s \in [0,r]}{\sup} \ \left\| x(r) \right\|^p \right] \; \mathrm{d}s \right) .
\end{align*}
Similar computations apply to the other terms and when considering solutions to \eqref{eq_SDE}. Therefore, we conclude from a routine Gr\"onwall inequality argument.

Let us prove the second inequality of the lemma. For $t \in [0,T]$, we compute
\begin{align*}
    \mathbb{E}\bigg[ \underset{s \in [0,t]}{\sup} \ &\| x_{u_1}(s) - x_{u_2}(s) \|^p \bigg] \le C \mathbb{E}\left[ \underset{s \in [0,t]}{\sup} \ \left\| \int^r_0 \left( b_0(s,x_{u_1}(s)) - b_0(s,x_{u_2}(s)) \right) \; \mathrm{d}s \right\|^p \right] \\
    & + C \sum^m_{i=1} \mathbb{E}\left[ \underset{s \in [0,t]}{\sup} \ \left\| \int^r_0 \left( u^i_1(s) b_i(s,x_{u_1}(s)) - u^i_2(s) b_i(s,x_{u_2}(s)) \right) \; \mathrm{d}s \right\|^p \right] \\
    & + C \mathbb{E}\left[ \underset{s \in [0,t]}{\sup} \ \left\| \int^r_0 \left( \sigma(s,x_{u_1}(s)) - \sigma(s,x_{u_2}(s)) \right) \; \mathrm{d}B_s \right\|^p \right] .
\end{align*}
For the second term on the right-hand side, for $i=1,\dots,m$ H\"older inequality gives
\begin{align*}
    &\mathbb{E}\left[ \underset{s \in [0,t]}{\sup} \ \left\| \int^r_0 \left( u^i_1(s) b_i(s,x_{u_1}(s)) - u^i_2(s) b_i(s,x_{u_2}(s)) \right) \; \mathrm{d}s \right\|^p \right] \le \\
    & \le C \left( \int^t_0 \mathbb{E}\left[ \underset{r \in [0,s]}{\sup} \ \left\| x_{u_1}(r) - x_{u_2}(r) \right\|^p \right] \; \mathrm{d}s + \mathbb{E}\left[ \left( \int^T_0 \| u_1(s) - u_2(s) \| \; \mathrm{d}s \right)^p \right] \right) ,
\end{align*}
and similar computations hold for the remaining terms and when considering solutions to \eqref{eq_SDE}. Again, we conclude by a Gr\"onwall inequality argument.
\end{proof}

\subsection{Proof of Lemma \ref{lemma_needle} and of Lemma \ref{lemma_needleGeneralized}}

The proof of Lemma \ref{lemma_needle} immediately follows from the following preliminary result.
\begin{lmm}[Stochastic needle-like variation formula - no free final time] \label{lemma_needlePrel}
Let $(\eta_1,\dots,\eta_j) \in \textnormal{Pr}_{\mathbb{R}^j_+}(\mathcal{C}^j_t)$ (in particular, no $(A'_2)$ or any assumption on $t > t_j$ are required). For $\varepsilon \in [0,t - t_j)$, uniformly for $\delta \in [-\varepsilon,T - t]$,
{\small
$$
\mathbb{E}\bigg[ \bigg\| \tilde g(\tilde x_{u_{\pi}}(t + \delta)) - \tilde g(\tilde x_u(t + \delta)) - \sum_{i=1}^j \eta_i \frac{\partial \tilde g}{\partial \tilde x}(\tilde x_u(t + \delta)) z_{t_i,u_i}(t + \delta) \bigg\|^2 \bigg]^{\frac{1}{2}} = o\left( \sum_{i=1}^j \eta_i \right) .
$$}
\end{lmm}

\begin{proof}
We only consider the case $j = 1$, being the most general case $j > 1$ done by adopting a classical induction argument (see, e.g., \cite{agrachev2013}). We only need to prove that
\begin{equation} \label{eq:needle1}
    \mathbb{E}\left[ \left\| \tilde x_{u_{\pi}}(t + \delta) - \tilde x_u(t + \delta) - \eta_1 z_{t_1,u_1}(t + \delta) \right\|^2 \right] = o(\eta^2_1)
\end{equation}
uniformly for every $\delta \in [-\varepsilon,T - t]$. For this, first we remark that $\tilde x_{u_{\pi}}(t) = \tilde x_u(t)$ for $t \in [0,t_1]$. Since $t \in (t_j,t_f]$ and $\varepsilon \ge 0$ are fixed and we take the limit $\eta_1 \rightarrow 0$, we may assume that $\eta_1 + t_1 < t - \varepsilon$. Therefore, without loss of generality, we replace $t + \delta$ with $t$, assuming $t \ge t_1 + \eta_1$ uniformly. We have (below, $C \ge 0$ denotes a constant)
\begin{align*}
    \mathbb{E}&\left[ \left\| \tilde x_{u_{\pi}}(t) - \tilde x_u(t) - \eta_1 z_{t_1,u_1}(t) \right\|^2 \right] \le C \mathbb{E}\left[ \left\| \int^{t_1 + \eta_1}_{t_1} \eta_1 A(s) z_{t_1,u_1}(s) \mathrm{d}s \right\|^2 \right] \\
    & + C \mathbb{E}\left[ \left\| \int^t_{t_1 + \eta_1} \left( \tilde b(s,u_{\pi}(s),x_{u_{\pi}}(s)) - \tilde b(s,u(s),x_u(s)) - \eta_1 A(s) z_{t_1,u_1}(s) \right) \mathrm{d}s \right\|^2 \right] \\
    & + C \mathbb{E}\left[ \left\| \int^t_0 \mathbbm{1}_{(t_1,T]}(s) \left( \tilde \sigma(s,x_{u_{\pi}}(s)) - \tilde \sigma(s,x_u(s)) - \eta_1 D(s) z_{t_1,u_1}(s) \right) \mathrm{d}B_s \right\|^2 \right] \\
    & + C \mathbb{E}\left[ \left\| \int^{t_1 + \eta_1}_{t_1} \left( \tilde b(s,u_{\pi}(s),x_{u_{\pi}}(s)) - \tilde b(s,u(s),x_u(s)) \right) \mathrm{d}s - \eta_1 z_{t_1,u_1}(t_1) \right\|^2 \right] .
\end{align*}
Let us analyze those integrals separately. Starting with the last one, we have
\begin{align*}
    \mathbb{E}&\left[ \left\| \int^{t_1 + \eta_1}_{t_1} \left( \tilde b(s,u_{\pi}(s),x_{u_{\pi}}(s)) - \tilde b(s,u(s),x_u(s)) \right) \mathrm{d}s - \eta_1 z_{t_1,u_1}(t_1) \right\|^2 \right] \le \\
    & \le C \mathbb{E}\left[ \left\| \int^{t_1 + \eta_1}_{t_1} \left( \tilde b(s,u_1,x_{u_{\pi}}(s)) - \tilde b(s,u_1,x_u(s)) \right) \mathrm{d}s \right\|^2 \right] \\
    & \quad + C \mathbb{E}\left[ \left\| \int^{t_1 + \eta_1}_{t_1} \left( \tilde b(s,u_1,x_u(s)) - \tilde b(s,u(s),x_u(s)) \right) \mathrm{d}s - \eta_1 z_{t_1,u_1}(t_1) \right\|^2 \right] .
\end{align*}
Lemma \ref{lemma_bound} immediately gives that the first term on the right-hand side is $o(\eta^2_1)$. In addition, from H\"older inequality, it follows that
\begin{align*}
    \frac{1}{\eta^2_1} \mathbb{E}&\left[ \left\| \int^{t_1 + \eta_1}_{t_1} \left( \tilde b(s,u_1,x_u(s)) - \tilde b(s,u(s),x_u(s)) \right) \mathrm{d}s - \eta_1 z_{t_1,u_1}(t_1) \right\|^2 \right] \le \\
    & \le \frac{C}{\eta_1} \int^{t_1 + \eta_1}_{t_1} \mathbb{E}\left[ \left\| \tilde b(s,u_1,x_u(s)) - \tilde b(t_1,u_1,x_u(t_1)) \right\|^2 \right] \; \mathrm{d}s \\
    & \quad + \frac{C}{\eta_1} \int^{t_1 + \eta_1}_{t_1} \mathbb{E}\left[ \left\| \tilde b(s,u(s),x_u(s)) - \tilde b(t_1,u(t_1),x_u(t_1)) \right\|^2 \right] \; \mathrm{d}s .
\end{align*}
Since $t_1$ is a Lebesgue point for $u$, the two terms above go to zero as $\eta_1 \rightarrow 0$. Next, by the Burkholder–Davis–Gundy inequality and a Taylor development, we have
\begin{align*}
    &\mathbb{E}\left[ \left\| \int^t_0 \mathbbm{1}_{(t_1,T]}(s) \left( \tilde \sigma(s,x_{u_{\pi}}(s)) - \tilde \sigma(s,x_u(s)) - \eta_1 D(s) z_{t_1,u_1}(s) \right) \mathrm{d}B_s \right\|^2 \right] \\
    & \le C \mathbb{E}\left[ \int^t_{t_1} \left\| D(s) \left( \tilde x_{u_{\pi}}(s) - \tilde x_u(s) - \eta_1 z_{t_1,u_1}(s) \right) \right\|^2 \mathrm{d}s \right] \\
    & \quad + C \mathbb{E}\left[ \int^t_{t_1} \left\| \int^1_0 \theta \frac{\partial^2 \tilde \sigma}{\partial \tilde x^2}\left( s, \theta \tilde x_u(s) + (1 - \theta) (\tilde x_{u_{\pi}} - \tilde x_u)(s) \right) \left( \tilde x_{u_{\pi}} - \tilde x_u \right)^2(s) \; \mathrm{d}\theta \right\|^2 \mathrm{d}s \right] \\
    & = C \int^t_{t_1} \mathbb{E}\left[ \left\| \tilde x_{u_{\pi}}(s) - \tilde x_u(s) - \eta_1 z_{t_1,u_1}(s) \right\|^2 \right] \; \mathrm{d}s + o(\eta^2_1) ,
\end{align*}
from Lemma \ref{lemma_bound}. Similar estimates hold for the remaining terms in the first inequality of the proof. Summarizing, there exists a constant $C \ge 0$ such that
{
$$
\mathbb{E}\Big[ \| \tilde x_{u_{\pi}}(t) - \tilde x_u(t) - \eta_1 z_{t_1,u_1}(t)  \|^2 \Big] \le o(\eta^2_1) + C \int^t_{t_1} \mathbb{E}\left[ \left\| \tilde x_{u_{\pi}}(s) - \tilde x_u(s) - \eta_1 z_{t_1,u_1}(s) \right\|^2 \right] \; \mathrm{d}s
$$}
and the conclusion follows from a routine Gr\"onwall inequality argument.
\end{proof}

\begin{proof}[Proof of Lemma \ref{lemma_needleGeneralized}]
We may assume $\delta \neq 0$ and $(A'_2)$. Hence, $\tilde g$ is an affine function, and in the rest of this proof we denote $\tilde g(\tilde x) = M \tilde x + d$, where $M \in \mathbb{R}^{(q+1) \times (n+1)}$ and $d \in \mathbb{R}^{q+1}$.

Developing, we have
\begin{align*}
    \bigg\| \mathbb{E}\bigg[ &\tilde g(\tilde x_{u_{\pi}}(t + \delta)) - \tilde g(\tilde x_u(t)) - \delta M \tilde b(t,u(t),x_u(t)) - \sum^j_{i=1} \eta_i A z_{t_i,u_i}(t) \bigg] \bigg\| \le \\
    & \le \Big\| \mathbb{E}\Big[ \tilde g(\tilde x_u(t + \delta)) - \tilde g(\tilde x_u(t)) - \delta M \tilde b(t,u(t),x_u(t)) \Big] \Big\| \\
    & \quad + \bigg\| \mathbb{E}\bigg[ \tilde g(\tilde x_{u_{\pi}}(t + \delta)) - \tilde g(\tilde x_u(t+\delta)) - \sum^j_{i=1} \eta_i M z_{t_i,u_i}(t+\delta) \bigg] \bigg\| + \sum^j_{i=1} \eta_i \left\| M \right\| \Big\| \mathbb{E}\Big[ z_{t_i,u_i}(t+\delta) - z_{t_i,u_i}(t) \Big] \Big\| .
\end{align*}
From Lemma \ref{lemma_needlePrel}, the second term on the right-hand side is $o\left( \sum_{i=1}^j \eta_i \right)$. For the last summand, from the property of the stochastic integral, we have
\begin{align*}
    \Big\| &\mathbb{E}\Big[ z_{t_i,u_i}(t+\delta) - z_{t_i,u_i}(t) \Big] \Big\| \le \left\| \mathbb{E}\left[ \int^{t + \delta}_t A(s) z_{t_i,u_i}(s) \; \mathrm{d}s \right] \right\| \\
    & + \left\| \mathbb{E}\left[ \int^{t + \delta}_0 \mathbbm{1}_{[t,T]}(s) D(s) z_{t_i,u_i}(s) \; \mathrm{d}B_s \right] \right\| \le C \mathbb{E}\left[ \underset{s \in [0,T]}{\sup} \ \| z_{t_i,u_i}(s) \| \; \right] \delta ,
\end{align*}
where the last term is $o\left(\delta + \sum_{i=1}^j \eta_i \right)$. It is worth pointing out the importance for $g$ being affine to provide the last claim. Finally, for the remaining term, we apply It\^o formula to each coordinate $h=1,\dots,q+1$, obtaining
{\small
\begin{align*}
    \Big[ \tilde g(\tilde x_u(t + \delta)) - \tilde g(\tilde x_u(t)) - &\delta M \tilde b(t,u(t),x_u(t)) \Big]_h = \sum^{n+1}_{k=1} M_{hk} \int^{t+\delta}_0 \mathbbm{1}_{[t,T]}(s) \tilde \sigma_k(s,x_u(s)) \; \mathrm{d}B_s \\
    & + \sum^{n+1}_{k=1} M_{hk} \left( \int^{t+\delta}_t \tilde b_k(s,u(s),x_u(s)) \; \mathrm{d}s - \delta \tilde b_k(t,u(t),x_u(t)) \right) ,
\end{align*}}
and therefore
\begin{align*}
    \frac{1}{\delta + \sum_{i=1}^j \eta_i} \Big\| &\mathbb{E}\Big[ \tilde g(\tilde x_u(t + \delta)) - \tilde g(\tilde x_u(t)) - \delta M \tilde b(t,u(t),x_u(t)) \Big] \Big\| \le \\
    & \le \frac{\| M \|}{\delta} \int^{t+\delta}_t \mathbb{E}\left[ \left\| \tilde b_k(s,u(s),x_u(s)) - \tilde b_k(t,u(t),x_u(t)) \right\| \right] \; \mathrm{d}s .
\end{align*}
As $t$ is Lebesgue for $u$, this last quantity goes to zero for $\delta \rightarrow 0$, and we conclude.
\end{proof}
\end{document}